\documentclass[12pt]{article}
\usepackage{amsmath,amssymb,amsthm}
\usepackage{graphicx}

\newtheorem{theorem}{Theorem}[subsection]
\newtheorem{proposition}[theorem]{Proposition}
\newtheorem{corollary}[theorem]{Corollary}
\newtheorem{lemma}[theorem]{Lemma}

\theoremstyle{definition}

\newtheorem*{corollary*}{Corollary}
\newtheorem{definition}[theorem]{Definition}
\newtheorem{example}[theorem]{Example}
\newtheorem{remark}[theorem]{Remark}
\newtheorem{terminology}[theorem]{Terminology}
\newtheorem{notation}[theorem]{Notation}

\newcommand{\pair}[2]{$\left[\scriptstyle\begin{matrix}#1\\#2\end{matrix}\right]$}
\DeclareMathOperator{\PL}{PL}
\DeclareMathOperator{\LL}{\mathbf{L}}
\DeclareMathOperator{\NN}{\mathbf{N}}
\DeclareMathOperator{\RR}{\mathbf{R}}
\DeclareMathOperator{\II}{\mathbf{I}}

\begin{document}
\title{Forest Diagrams for Elements of \\Thompson's Group $F$}
\author{James M. Belk and Kenneth S. Brown}
\maketitle

\begin{abstract}
We introduce \emph{forest diagrams} to represent elements of Thompson's group $F$. These diagrams relate
to a certain action of $F$ on the real line in the same way that tree diagrams relate to the standard
action of $F$ on the unit interval. Using forest diagrams, we give a conceptually simple length formula
for elements of $F$ with respect to the $\{x_0,x_1\}$ generating set, and we discuss the construction
of minimum-length words for positive elements.  Finally, we use forest diagrams and the length formula
to examine the structure of the Cayley graph of $F$.
\end{abstract}

\section{Introduction}

Thompson's group $F$ is defined by the following infinite presentation:
\begin{equation*}
F=\langle x_0,x_1,x_2,\ldots\mid x_nx_k=x_kx_{n+1}\,\text{for }n>k\rangle
\end{equation*}
It is isomorphic to the group $\PL_2\left(I\right)$ of all piecewise-linear,
orientation-preserving homeomorphisms of the unit interval satisfying the
following conditions:
\begin{enumerate}
\item All slopes are integral powers of $2$, and
\item All breakpoints have dyadic rational coordinates.
\end{enumerate}
The group $F$ was first studied by Richard J. Thompson in the 1960s. The
standard introduction to $F$ is \cite{CFP}.

This paper is organized as follows:

\begin{itemize}
\item In \textbf{Section 2}, we give the necessary background regarding $F$. In
particular, we review how elements of $\PL_2(I)$ can be described
by tree diagrams.

\item In \textbf{Section 3}, we introduce a group $\PL_2(\mathbb{R})$ of
piecewise-linear homeomorphisms of the real line that is isomorphic
with $F$. We then show how to represent elements of $\PL(\mathbb{R})$ by
\emph{forest diagrams}.

\item In \textbf{Section 4}, we use forest diagrams to examine the lengths of elements
of Thompson's group with respect to the $\{x_0,x_1\}$ generating
set. We begin by studying positive elements, where the situation is quite
simple, and then move on to the general length formula.

\item In \textbf{Section 5}, we give some further applications of forest diagrams
and the length formula.
\end{itemize}

\

\section{Background on $F$}

Most of the results in this section are stated without proof. Details
can be found in \cite{CFP}.

\

\subsection{Tree Diagrams}
Suppose we take the interval $[0,1]$ and cut it in half, like this:
\begin{center}
\includegraphics{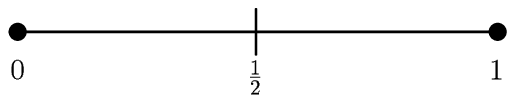}
\end{center}
We then cut each of the resulting intervals in half:
\begin{center}
\includegraphics{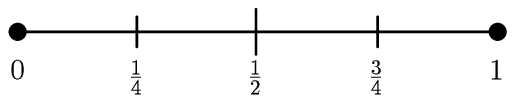}
\end{center}
and then cut some of the new intervals in half:
\begin{center}
\includegraphics{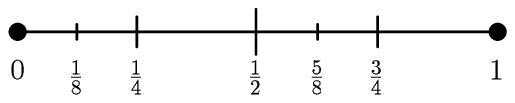}
\end{center}
to get a certain subdivision of $[0,1]$. Any subdivision of $[0,1]$
obtained in this manner (i.e. by repeatedly cutting intervals in half)
is called a \emph{dyadic subdivision}.

The intervals of a dyadic subdivision are all of the form:
\begin{equation*}
\left[\frac{k}{2^n},\frac{k+1}{2^n}\right]\qquad k,n\in\mathbb{N}
\end{equation*}
These are the \emph{standard dyadic intervals}. We could alternatively define a
dyadic subdivision as any partition of $[0,1]$ into standard dyadic intervals.

Each element of $\PL_2(I)$ can be described by a pair of dyadic subdivisions:

\begin{proposition}
Let $f\in\PL_2(I)$. Then there exist dyadic subdivisions
$\mathcal{D},\mathcal{R}$ of $[0,1]$ such that $f$ maps each interval of
$\mathcal{D}$ linearly onto an interval of $\mathcal{R}$.\quad\qedsymbol
\end{proposition}

\begin{example}
Consider the element $f\in\PL_2(I)$ with graph:
\begin{center}
\includegraphics{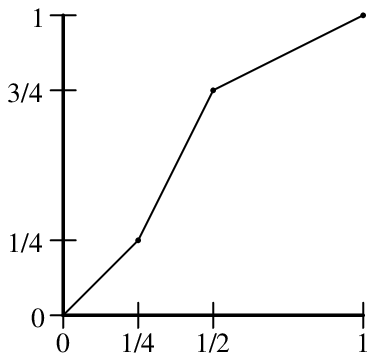}
\end{center}
Then $f$ maps intervals of the subdivision:
\begin{center}
\includegraphics{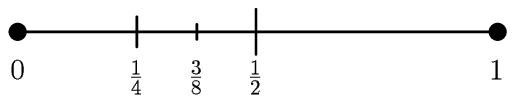}
\end{center}
linearly onto intervals of the subdivision:
\begin{center}
\includegraphics{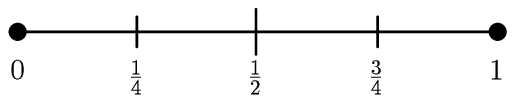}\\[18pt]
\end{center}
\end{example}

We can represent dyadic subdivisions of $[0,1]$ by finite binary trees.
For example, the subdivision:
\begin{center}
\includegraphics{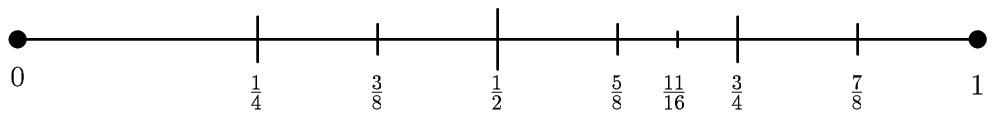}
\end{center}
corresponds to the binary tree:
\begin{center}
\includegraphics{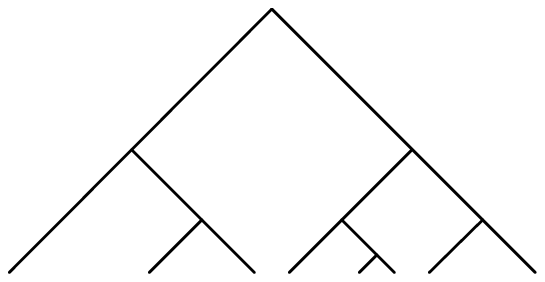}
\end{center}
Each leaf of this tree represents an interval of the subdivision, and the
root represents the interval $[0,1]$. The other nodes represent standard
dyadic intervals from intermediate stages of the dyadic subdivision.

Combining this observation with proposition 2.1.1, we see that any element
$f\in\PL_2(I)$ can be described by a pair of binary trees. This is called
a \emph{tree diagram} for $f$.

\begin{example}
Let $f$ be the element of $\PL_2(I)$ from example 2.1.2. Then $f$ has tree diagram:
\begin{center}
\includegraphics{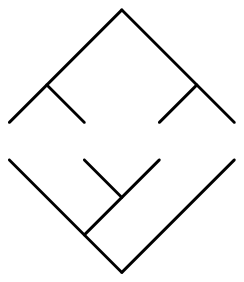}
\end{center}
We have aligned the two trees vertically so that corresponding leaves match up.
By convention, the \emph{domain tree} appears on the \emph{bottom}, and the
\emph{range tree} appears on the \emph{top}.
\end{example}

The tree diagram for an element $f\in\PL_2(I)$ is not unique. For example,
all of the following are tree diagrams for the identity:
\begin{center}
\includegraphics{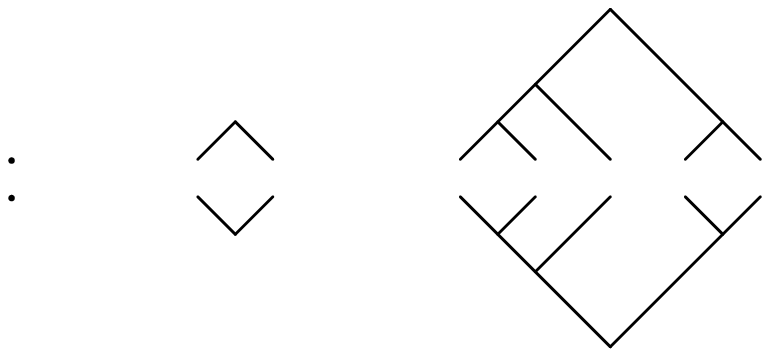}
\end{center}
In general, a \emph{reduction} of a tree diagram consists of removing an opposing
pair of carets, like this:
\begin{center}
\includegraphics{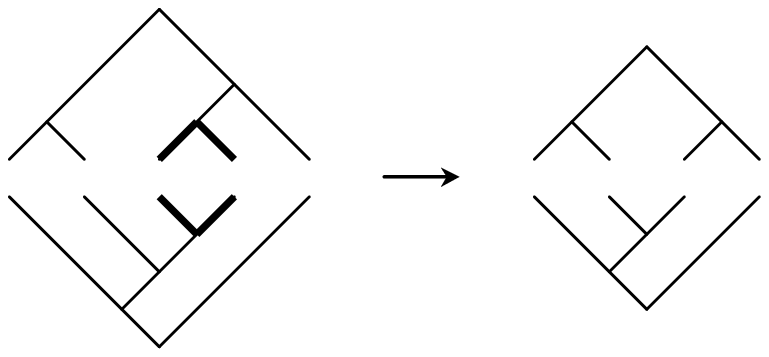}
\end{center}
Performing a reduction does not change the element of $\PL_2(I)$
described by a tree diagram: it merely corresponds to removing an
unnecessary ``cut'' from the subdivisions of the domain and range.

\begin{definition}
A tree diagram is \emph{reduced} if it has no opposing pairs of carets.
\end{definition}

\begin{proposition}
Every element of $\PL_2(I)$ has a unique reduced tree diagram.\quad\qedsymbol
\end{proposition}

\

\subsection{Positive Elements and Normal Form}

Recall that $F$ has presentation:
\begin{equation*}
F=\langle x_0,x_1,x_2,\ldots\mid x_nx_k=x_kx_{n+1}\text{ for }k<n\rangle
\end{equation*}
We have previously asserted that $F$ is isomorphic with $\PL_2(I)$. One
such isomorphism is defined as follows:
\begin{center}
\includegraphics{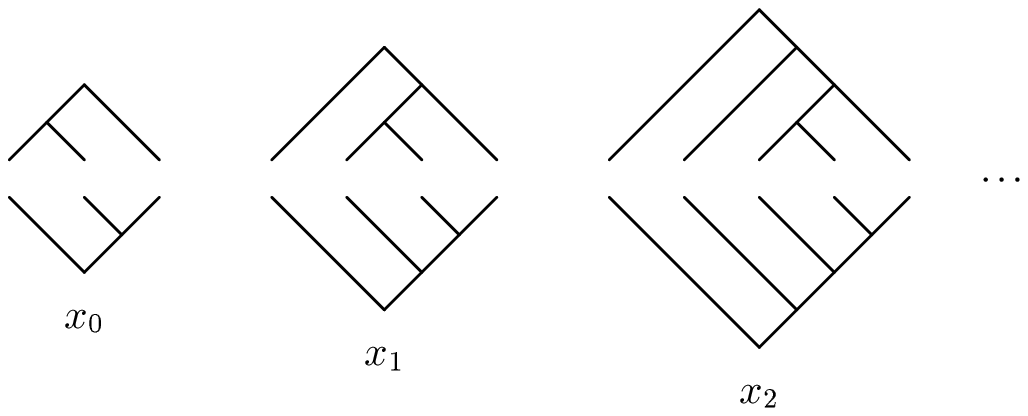}
\end{center}
Note that the domain trees of the $x_i$'s all have the property that no caret
has a left child. Such a tree is called a \emph{right vine}:
\begin{center}
\includegraphics{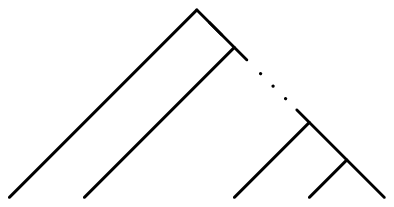}
\end{center}

\begin{definition}
An element of $F$ is \emph{positive} if it lies in the submonoid generated by
$\{x_0,x_1,x_2,\ldots\}$.
\end{definition}

\begin{proposition}
An element of $F$ is positive if and only if the bottom tree of its reduced tree
diagram is a right vine.\quad\qedsymbol
\end{proposition}

It turns out that $F$ is the group of fractions of its positive monoid, in the
sense that any element of $F$ can be written as $pq^{-1}$ for some positive $p$~and~$q$.
More precisely:

\begin{proposition}[Normal Form]
Every element of $F$ can be expressed uniquely in the form:
\begin{equation*}
x_0^{a_0}\cdots x_n^{a_n}x_n^{-b_n}\cdots x_0^{-b_0}
\end{equation*}
where $a_0,\ldots,a_n,b_0,\ldots,b_n\in\mathbb{N}$ and:
\begin{enumerate}
\item Either $a_n>0$ or $b_n>0$, but not both.
\item If both $a_i>0$ and $b_i>0$, then either $a_{i+1}>0$ or $b_{i+1}>0$.
\quad\qedsymbol\renewcommand{\qedsymbol}{}
\end{enumerate}
\end{proposition}

The first half of the normal form is called the \emph{positive part} of an element,
and the second half is called the \emph{negative part}. These halves correspond to
the two halves of the tree diagram:

\begin{proposition}
Let \pair{T_+}{T_-} be the reduced tree diagram for an element $f\in F$, and let $V$
be a right vine with the same number of leaves as $T_+$ and $T_-$. Then \pair{T_+}{V}
is a tree diagram for the positive part of $f$, and \pair{V}{T_-} is a tree diagram
for the negative part of $f.\quad\qedsymbol$
\end{proposition}

\

\section{Forest Diagrams}

It is immediate from the presentation of $F$ that:
\begin{equation*}
x_n=x_0^{1-n}x_1x_0^{n-1}
\end{equation*}
for all $n\geq 1$. \ Therefore, $F$ is generated by the two elements $\{x_0,x_1\}$.

In this section, we describe a group $\PL_2(\mathbb{R})$ of self-homeomorphisms of the real line
that is isomorphic to $F$, and develop \emph{forest diagrams} in analogy with the
development of tree diagrams in the previous section.  These forest diagrams seem to interact
particularly nicely with the $\{x_0,x_1\}$-generating set.

The existence of forest diagrams was noted by K. Brown in \cite{Bro},
but the pictures themselves have not previously appeared in the literature.
They are similar to the ``diagrams'' of V. Guba and M. Sapir (see \cite{GuSa} and \cite{Guba}).

\

\subsection{The Group $\PL_2(\mathbb{R})$}

Let $\PL_2(\mathbb{R})$ be the group of all piecewise-linear, orientation-preserving
self-homeomorphisms $f$ of $\mathbb{R}$ satisfying the following conditions:
\begin{enumerate}
\item Each linear segment of $f$ has slope a power of $2$.
\item $f$ has only finitely many breakpoints, each of which has dyadic rational coordinates.
\item The leftmost linear segment of $f$ is of the form:
\begin{equation*}
f(t)=t-m
\end{equation*}
and the rightmost segment is of the form:
\begin{equation*}
f(t)=t-n
\end{equation*}
for some integers $m,n$.
\end{enumerate}
The following is well-known:

\begin{proposition}
$\PL_2(\mathbb{R})$ is isomorphic with $\PL_2(I)$.
\end{proposition}

\begin{proof}
Let $\psi\colon\mathbb{R}\rightarrow(0,1)$ be the
piecewise-linear homeomorphism that maps the intervals:
\begin{center}
\includegraphics{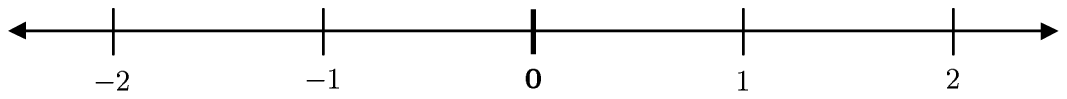}
\end{center}
linearly onto the intervals:
\begin{center}
\includegraphics{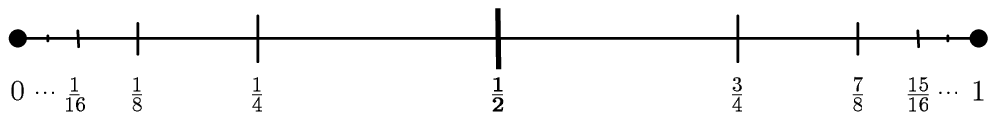}
\end{center}
Then $f\mapsto\psi^{-1}f\psi$ is the desired isomorphism $\PL_2(I)\rightarrow\PL_2(\mathbb{R})$.
\end{proof}

\begin{corollary}
$\PL_2(\mathbb{R})$ is isomorphic with $F$. The generators $\{x_0,x_1\}$ of $F$ map to the functions:
\begin{equation*}
x_0(t)=t-1
\end{equation*}
and:
\begin{center}
$x_1(t)=
\begin{cases}
t&t\leq 0\\
\frac{1}{2}t&0\leq t\leq 2\\
t-1&t\geq 2
\end{cases}
\qquad\qquad$\parbox{1.5in}{\includegraphics{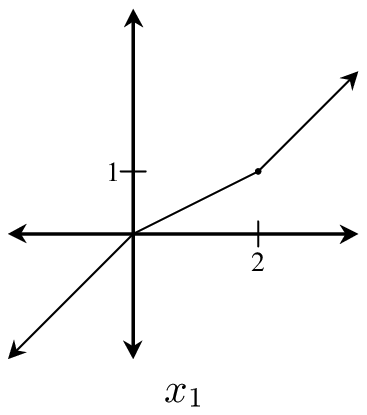}}
\end{center}
\end{corollary}

\

\subsection{Forest Diagrams for Elements of $\PL_2(\mathbb{R})$}

We think of the real line as being pre-subdivided as follows:
\begin{center}
\includegraphics{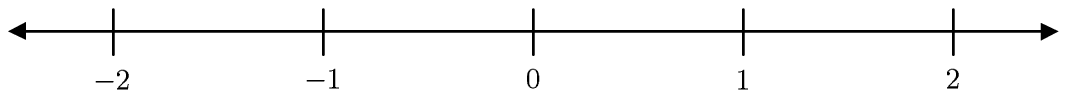}
\end{center}
A \emph{dyadic subdivision} of $\mathbb{R}$ is a subdivision obtained by cutting
finitely many of these intervals in half, and then cutting finitely many of the resulting
intervals in half, etc.

\begin{proposition}
Let $f\in\PL_2(\mathbb{R})$. Then there exist dyadic subdivisions
$\mathcal{D},\mathcal{R}$ of $\mathbb{R}$ such that $f$ maps each interval of
$\mathcal{D}$ linearly onto an interval of $\mathcal{R}.\quad\qedsymbol$
\end{proposition}

A \emph{binary forest} is a sequence $(\ldots,T_{-1},T_0,T_1,\ldots)$ of
finite binary trees. We depict such a forest as a line of binary trees
together with a pointer at~$T_0$:
\begin{center}
\includegraphics{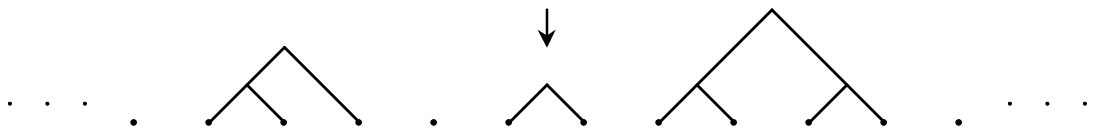}
\end{center}
A binary forest is \emph{bounded} if only finitely many of the trees
$T_i$ are nontrivial.

Every bounded binary forest corresponds to some dyadic subdivision of the real line.
For example, the forest above corresponds to the subdivision:
\begin{center}
\includegraphics{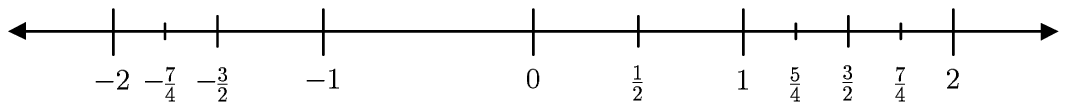}
\end{center}
Each tree $T_i$ represents an interval $[i,i+1]$, and each leaf represents an
interval of the subdivision.

Combining this with proposition 3.2.1, we see that any $f\in\PL_2(\mathbb{R})$
can be represented by a pair of bounded binary forests, together with an
order-preserving bijection of their leaves. This is called a \emph{forest diagram}
for $f$.

\begin{example}
Let $f$ be the element of $\PL_2(\mathbb{R})$ with graph:
\begin{center}
\includegraphics{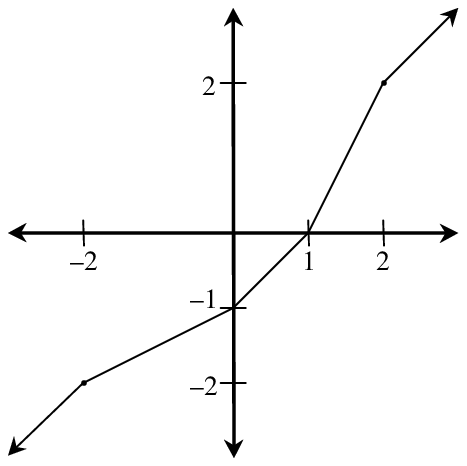}
\end{center}
Then $f$ has forest diagram:
\begin{center}
\includegraphics{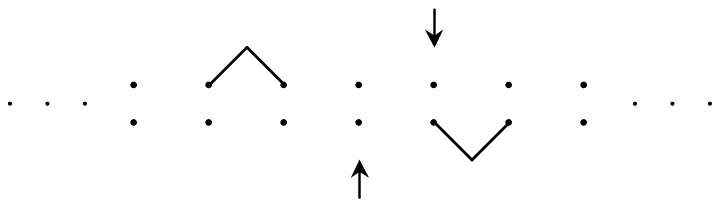}
\end{center}
Again, we have aligned the two forests vertically so that corresponding leaves match up.
By convention, the \emph{domain} tree appears on the \emph{bottom}, and the \emph{range}
tree appears on the \emph{top}.
\end{example}

\begin{example}
Here are the forest diagrams for $x_0$ and $x_1$:
\begin{center}
$x_0$:\qquad\parbox{2.5in}{\fbox{\includegraphics[clip]{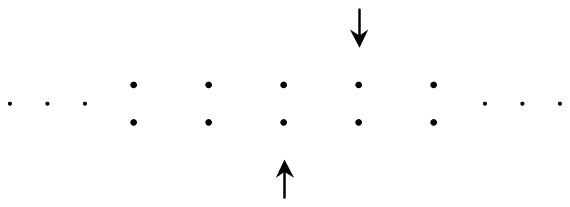}}}

\

\

$x_1$:\qquad\parbox{2.5in}{\fbox{\includegraphics[clip]{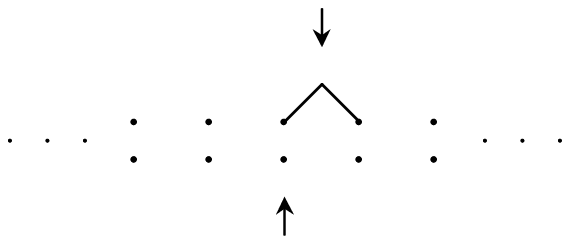}}}\\[18pt]
\end{center}
\end{example}

Of course, there are several forest diagrams for each element of $\PL_2(\mathbb{R})$.
In particular, it is possible to delete an opposing pair of carets:
\begin{center}
\includegraphics{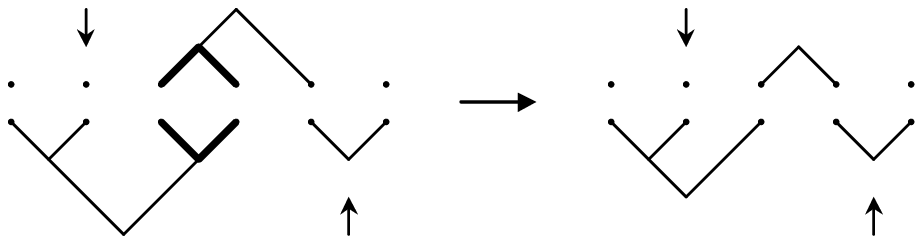}
\end{center}
without changing the resulting homeomorphism. This is called a \emph{reduction} of a forest diagram.
A forest diagram is \emph{reduced} if it does not have any opposing pairs of carets.

\begin{proposition}
Every element of $\PL_2(\mathbb{R})$ has a unique reduced forest diagram.\quad\qedsymbol
\end{proposition}

\begin{remark}
From this point forward, we will only draw the \emph{support} of the
forest diagram (i.e. the minimum interval containing both pointers
and all nontrivial trees), and we will omit the ``$\cdots$'' indicators.
\end{remark}

\begin{remark}
It is fairly easy to translate between tree diagrams and forest diagrams. Given a tree diagram:
\begin{center}
\includegraphics{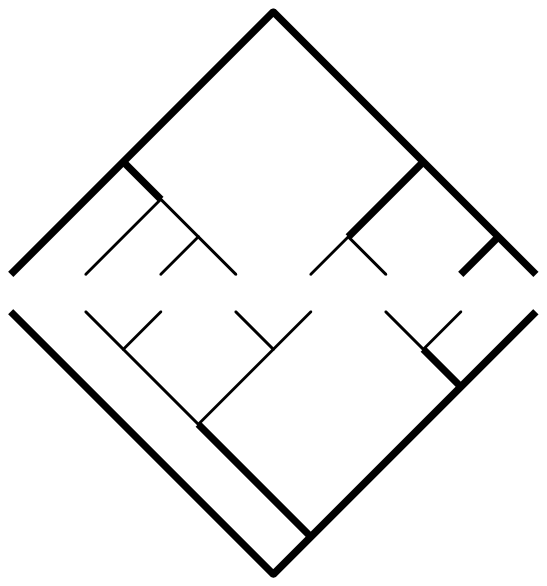}
\end{center}
we simply remove the outer layer of each tree to get the corresponding forest diagram:
\begin{center}
\includegraphics{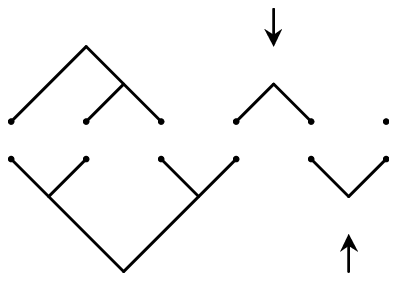}
\end{center}
Notice how the positions of the pointers in the forest diagram are determined by the former
positions of the roots in the tree diagram.
\end{remark}

\

\subsection{The Action of $\{x_0,x_1\}$}

The action of $\{x_0,x_1\}$ on forest diagrams is particularly nice:

\begin{proposition}
Let $\mathfrak{f}$ be a forest diagram for some $f\in F$. Then:
\begin{enumerate}
\item A forest diagram for $x_0f$ can be obtained by moving the
top pointer of $\mathfrak{f}$ one tree to the right.
\item A forest diagram for $x_1f$ can be obtained by attaching a
caret to the roots of the $0$-tree and $1$-tree in the top forest of
$\mathfrak{f}$. Afterwards, the top pointer points to the new, combined tree.\quad\qedsymbol
\end{enumerate}
\end{proposition}

If $\mathfrak{f}$ is reduced, then the given forest diagram for $x_0f$ will
always be reduced. The forest diagram given for $x_1f$ will not be reduced,
however, if the caret that was created opposes a caret from the bottom tree.
In this case, left-multiplication by $x_1$ effectively ``cancels'' the
bottom caret.

\begin{example}
Let $f\in F$ have forest diagram:
\begin{center}
\includegraphics{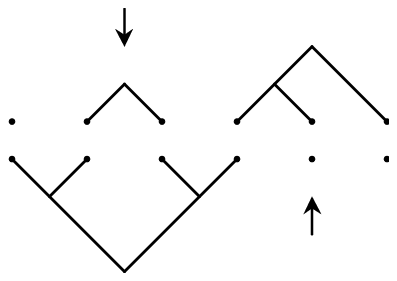}
\end{center}
Then $x_0f$ has forest diagram:
\begin{center}
\includegraphics{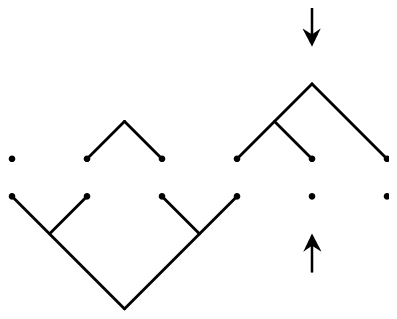}
\end{center}
and $x_1f$ has forest diagram:
\begin{center}
\includegraphics{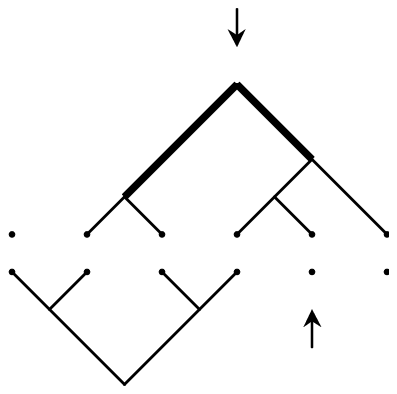}\\[18pt]
\end{center}
\end{example}

\pagebreak[0]
\begin{example}
Let $f\in F$ have forest diagram:
\begin{center}
\includegraphics{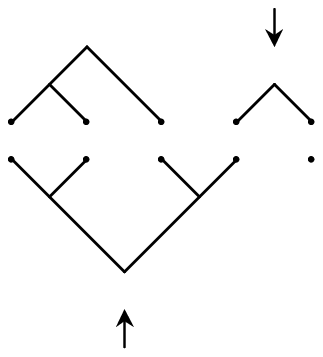}
\end{center}
Then $x_0f$ has forest diagram:
\begin{center}
\includegraphics{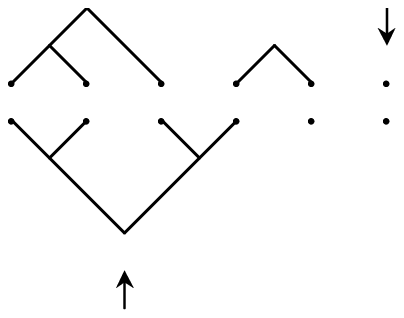}
\end{center}
and $x_1f$ has forest diagram:
\begin{center}
\includegraphics{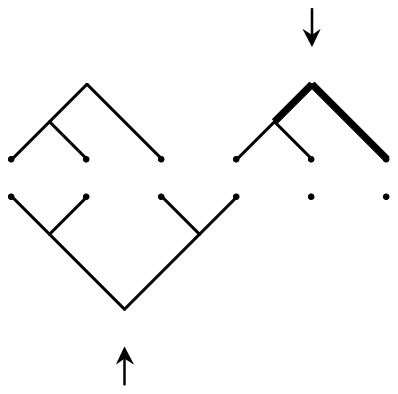}
\end{center}
Note that the forest diagrams for $x_0f$ and $x_1f$ both have larger support
than the forest diagram for $f$.
\end{example}

\begin{example}
Let $f\in F$ have forest diagram:
\begin{center}
\includegraphics{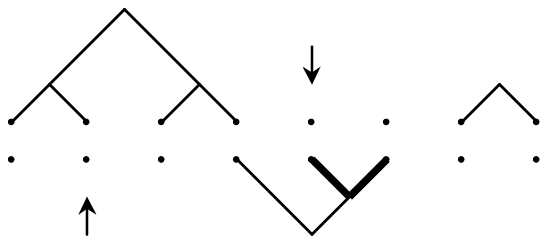}
\end{center}
Then $x_1f$ has forest diagram:
\begin{center}
\includegraphics{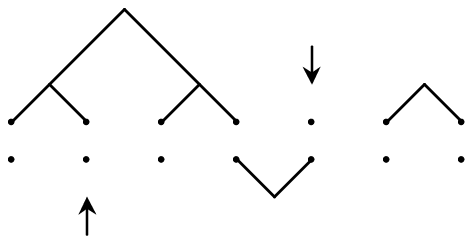}
\end{center}
Note that left-multiplication by $x_1$ canceled the highlighted bottom caret.
\end{example}

\begin{proposition}
Let $\mathfrak{f}$ be a forest diagram for some $f\in F$. Then:
\begin{enumerate}
\item A forest diagram for $x_0^{-1}f$ can be obtained by moving the
top pointer of $\mathfrak{f}$ one tree to the left.
\item A forest diagram for $x_1^{-1}f$ can be obtained by ``dropping a
negative caret'' at the current position of the top pointer. If the current tree
is nontrivial, the negative caret cancels with the top caret of the current tree,
and the pointer moves to the resulting left child. If the current tree is trivial,
the negative caret ``falls through'' to the bottom forest, attaching to the specified leaf.
\quad\qedsymbol
\end{enumerate}
\end{proposition}

\begin{example}
\quad Let $f$ and $g$ be the elements of $F$ with forest diagrams:
\begin{center}
\includegraphics{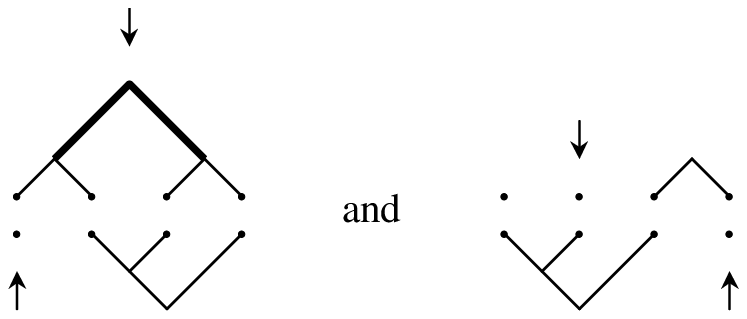}
\end{center}
Then $x_1^{-1}f$ and $x_1^{-1}g$ have forest diagrams:
\begin{center}
\includegraphics{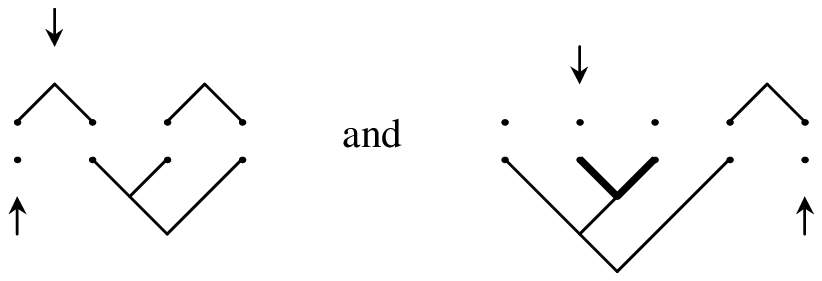}
\end{center}
In the first case, the $x_1^{-1}$ simply removed a caret from the top tree.
In the second case, there was no caret on top to remove, so a new caret was
attached to the leaf on the bottom. Note that this creates a new column
immediately to the right of the pointer.
\end{example}

\

\subsection{Positive Elements and Normal Form}

There is a close relationship between the normal form of an element and its forest diagram.
It hinges on the following proposition:

\begin{proposition}
Let $\mathfrak{f}$ be the forest diagram for some $f\in F$, and let $n>1$.  Then a forest
diagram for $x_nf$ can be obtained by attaching a caret to the roots of $T_{n-1}$ and $T_n$
in the top forest of $\mathfrak{f}$.
\end{proposition}
\begin{proof}
For $n>1$, $x_n = x_0^{1-n} x_1 x_0^{n-1}$.
\end{proof}

\begin{corollary}
Let $f\in F$, and let $\mathfrak{f}$ be its reduced forest diagram.  Then $f$ is positive
if and only if:
\begin{enumerate}
\item The entire bottom forest of $\mathfrak{f}$ is trivial, and
\item The bottom pointer is at the left end of the support of $\mathfrak{f}$.
\end{enumerate}
\end{corollary}

Using proposition 3.4.1, it is easy to construct the forest diagram for any positive element.
It is also possible to find the normal form when given the forest diagram:

\begin{example}
Suppose $f\in F$ has forest diagram:
\begin{center}
\includegraphics{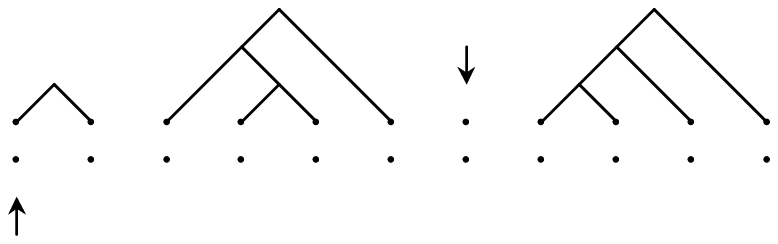}
\end{center}
Then:
\begin{equation*}
f = x_0^2 x_1 x_3^2 x_4 x_8^3
\end{equation*}
Since the top pointer of $f$ is two trees from the left, the normal form of $f$ has an $x_0^2$.
The powers of the other generators are determined by the number of carets built upon the
corresponding leaf.  Note that the carets are constructed from right to left.
\end{example}

It is not much harder to deal with mixed (non-positive) elements:

\pagebreak[2]
\begin{example}
The element:
\begin{equation*}
x_0^3 x_2 x_5^2 x_7 x_6^{-1} x_5^{-1} x_1^{-2} x_0^{-1}
\end{equation*}
has forest diagram:
\begin{center}
\includegraphics{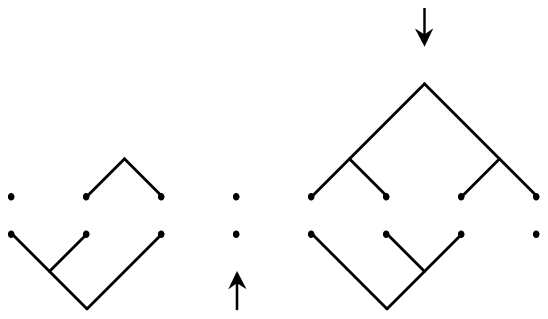}
\end{center}
\end{example}

\

\section{Lengths in $F$}

In this section, we derive a formula for the lengths of elements of $F$ with
respect to the $\{x_0,x_1\}$-generating set. This formula uses the forest
diagrams introduced in section 3.

Lengths in $F$ were first studied by S. B. Fordham. In his unpublished thesis \cite{Ford}
Fordham gives a formula for the length of an element of $F$ based on its tree diagram.  (See
\cite{ClTa1}, \cite{ClTa2}, or \cite{Bur} for a published version of Fordham's result.)  Our
length formula can be viewed as a simplification of Fordham's original work.

V. Guba has recently obtained another length formula for $F$ using the ``diagrams'' of Guba
and Sapir.  See \cite{Guba} for details.

\

\subsection{Lengths of Strongly Positive Elements}

We shall begin by investigating the lengths of strongly positive elements.
The goal is to develop some intuition for lengths before the statement of
the general length formula in section 4.2.

An element $f\in F$ is \emph{strongly positive} if it lies in the submonoid
generated by $\{x_1,x_2,\ldots\}$. Here is a forest diagram for a typical
strongly positive element:
\begin{center}
\includegraphics{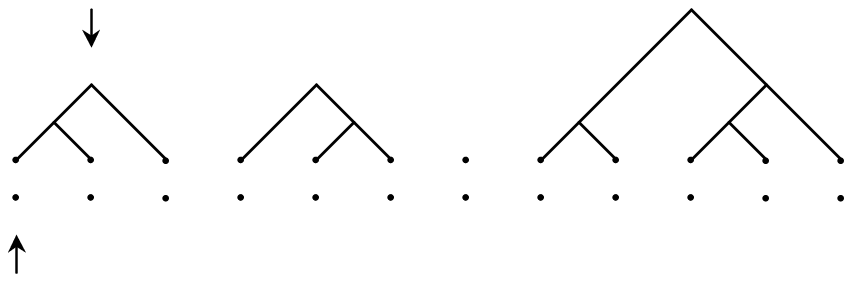}
\end{center}
Note that the entire bottom forest is trivial, and that both pointers are
at the left end of the support of $f$.

Logically, the results of this section depend on the general length formula.
In particular, we need the following lemma:

\begin{lemma}
Let $f\in F$ be strongly positive. Then there exists a minimum-length
word for $f$ with no appearances of $x_1^{-1}$.
\end{lemma}

This lemma is intuitively obvious: there should be no reason to ever create
bottom carets, or to delete top carets, when constructing a strongly positive
element.  Unfortunately, it would be rather tricky to supply a proof of this
fact.  Instead we refer the reader to corollary 4.3.8, from which the lemma follows
immediately.

From this lemma, we see that any strongly positive element $f\in F$ has a
minimum-length word of the form:
\begin{equation*}
x_0^{a_n}x_1\cdots x_0^{a_1}x_1x_0^{a_0}
\end{equation*}
where $a_0,\ldots,a_n\in\mathbb{Z}$. Since $f$ is strongly positive, we have: 
\begin{equation*}
a_0+\cdots+a_n=0
\end{equation*}
and
\begin{equation*}
a_0+\cdots+a_i\geq 0
\parbox{0in}{\parbox{2in}{\qquad(for $i=0,\ldots,n-1$)}}
\end{equation*}
Such words can be represented by words in $\{x_1,x_2,\ldots\}$ via
the identifications $x_n=x_0^{1-n}x_1x_0^{n-1}$. For example, the word:
\begin{equation*}
x_0^{-5}\,x_1\,x_0^{-2}\,x_1\,x_0^4\,x_1\,x_0^{-3}\,x_1\,x_0^6
\end{equation*}
can be represented by:
\begin{equation*}
x_6\,x_8\,x_4\,x_7
\end{equation*}
More generally:

\begin{notation}
We will use the word:
\begin{equation*}
x_{i_n}\cdots x_{i_2}x_{i_1}
\end{equation*}
in $\{x_1,x_2,\ldots\}$ to represent the word:
\begin{equation*}
x_0^{1-i_n}\,x_1\,\cdots\,x_0^{i_3-i_2}\,x_1\,x_0^{i_2-i_1}\,x_1\,x_0^{i_1-1}
\end{equation*}
in $\{x_0,x_1\}$.
\end{notation}

Note then that $x_{i_n}\cdots x_{i_2}x_{i_1}$ represents a word with length:
\begin{equation*}
\left(|1-i_n|+\cdots+|i_3-i_2|+|i_2-i_1|+|i_1-1|\right)+n
\end{equation*}
We now proceed to some examples, from which we will derive a general theorem.

\begin{example}
Let $f\in F$ be the element with forest diagram:
\begin{center}
\includegraphics{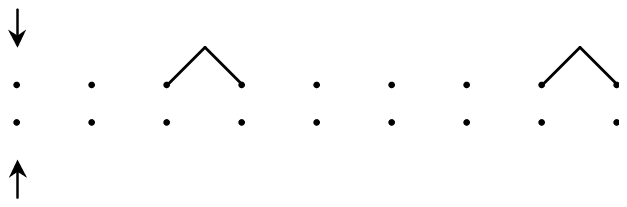}
\end{center}
There are only two candidate minimum-length words for $f$: $x_3x_8$ and
$x_7x_3$. Their lengths are:
\begin{equation*}
\begin{aligned}[t]
&&\quad&\left(2+5+7\right)+2=16&\quad&\quad\text{for the word }x_3x_8\\
&\text{and}\quad&\quad&\left(6+4+2\right)+2=14&\quad&\quad\text{for the word
}x_7x_3.
\end{aligned}
\end{equation*}
Let's see if we can explain this. The word $x_3x_8=x_0^{-2}x_1x_0^{-5}x_1x_0^7$
corresponds to the following construction of $f$:
\begin{enumerate}
\item Starting at the identity, move right seven times
and construct the right caret.
\item Next move left five times, and construct the left caret.
\item Finally, move left twice to position of the bottom pointer.
\end{enumerate}
This word makes a total of fourteen moves, crossing twice over each of seven spaces:
\begin{center}
\includegraphics{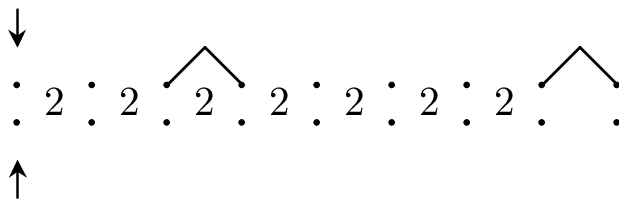}
\end{center}

On the other hand, the word $x_7x_3=x_0^{-6}x_1x_0^4x_1x_0^2$ corresponds to
the following construction:
\begin{enumerate}
\item Starting at the identity, move right twice and construct the
left caret.
\item Next move right four more times, and construct the right caret.
\item Finally, move left six times to the position of the bottom pointer.
\end{enumerate}
This word makes only twelve moves:
\begin{center}
\includegraphics{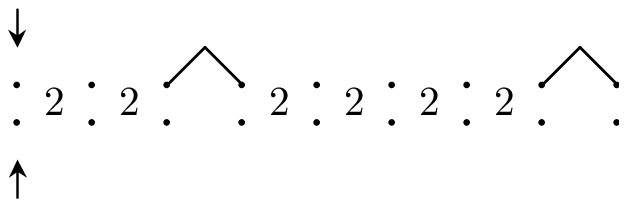}
\end{center}
In particular, this word never moves across the space under the left caret.
It avoids this by \emph{building the left caret early}.  Once the left caret
is built, the word can simply pass over the space under the left caret without
spending time to move across it.
\end{example}

\begin{terminology}
We call a space in a forest \emph{interior} if it lies under a tree (or
over a tree, if the forest is upside-down) and \emph{exterior} if it lies
between two trees.
\end{terminology}

\begin{example}
Let $f\in F$ be the element with forest diagram:
\begin{center}
\includegraphics{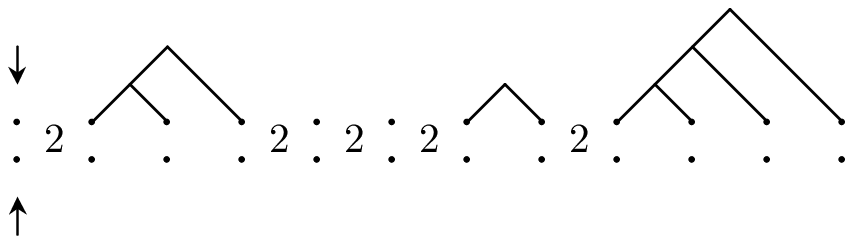}
\end{center}
Clearly, each of the five exterior spaces in the support of $f$ must be
crossed twice during construction. Furthermore, it is possible to avoid
crossing any of the interior spaces by \emph{constructing carets from left to
right}. In particular:
\begin{equation*}
x_6^3\,x_5\,x_2^2
\end{equation*}
is a minimum-length word for $f$. Therefore, $f$ has length:
\begin{equation*}
(5+1+3+1)+6=16\\[6pt]
\end{equation*}
\end{example}

It is not always possible to avoid crossing all the interior spaces:

\begin{example}
Let $f\in F$ be the element with forest diagram:
\begin{center}
\includegraphics{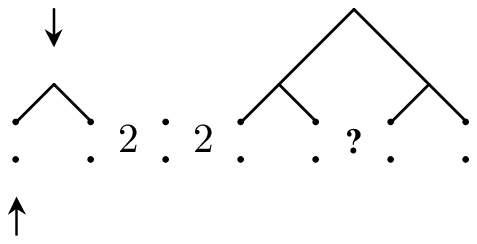}
\end{center}
Clearly, each of the two exterior spaces in the support of $f$ must be
crossed twice during construction. However, the space marked (?) must also
be crossed twice, since we must create the caret immediately to its right
before we can create the caret above it.

It turns out that these are the only spaces which must be crossed. For
example, the word:
\begin{equation*}
x_3\,x_4\,x_3\,x_1
\end{equation*}
crosses only these spaces.  Therefore, $f$ has length:
\begin{equation*}
(2+1+1+2+0)+4=10\\[6pt]
\end{equation*}
\end{example}

In this last example, we learned that it is not always possible to construct
carets from left to right. However, if one always constructs the \emph{leftmost
possible caret first}, then it is never necessary to move more than one space
to the left in the middle of the construction. This is the content of the
following theorem:

\begin{theorem}[Anti-Normal Form]
Let $f\in F$ be strongly positive. Then $f$ can be expressed uniquely
in the form:
\begin{equation*}
x_{i_n}\cdots x_{i_2}x_{i_1}
\end{equation*}
where $i_{k+1}\geq i_k-1$ for all $k.\quad\qedsymbol$
\end{theorem}

We say that a word:
\begin{equation*}
x_{i_n}\cdots x_{i_2}x_{i_1}
\end{equation*}
is in \emph{anti-normal form} if $i_{k+1}\geq i_k-1$ for each $k$. On the
forest diagram, anti-normal form corresponds to constructing the
\emph{leftmost possible caret} at each stage.

In contrast, the normal form for an element satisfies $i_{k+1}\leq i_k$ for
each $k$, and corresponds to constructing the rightmost possible caret at
each stage. This explains our terminology.

The anti-normal form for a strongly positive element of $F$ is clearly
minimum-length, since it crosses only those spaces in the forest diagram
that must be crossed. We can give an explicit length formula by counting
these spaces:

\begin{theorem}
Let $f\in F$ be strongly positive. Then the length of $f$ is:
\begin{equation*}
2\,n(f)+c(f)
\end{equation*}
where $n(f)$ and $c(f)$ are defined as follows. Let $\mathfrak{f}$ be the reduced
forest diagram for $f$. Then:
\begin{enumerate}
\item $n\left(f\right)$ is the number of spaces in the support of $\mathfrak{f}$
that are either exterior or lie immediately to the left of some caret, and
\item $c\left(f\right)$ is the number of carets in $\mathfrak{f}$.\quad\qedsymbol
\end{enumerate}
\end{theorem}

\begin{example}
\quad Let $f\in F$ be the element with forest diagram:
\begin{center}
\includegraphics{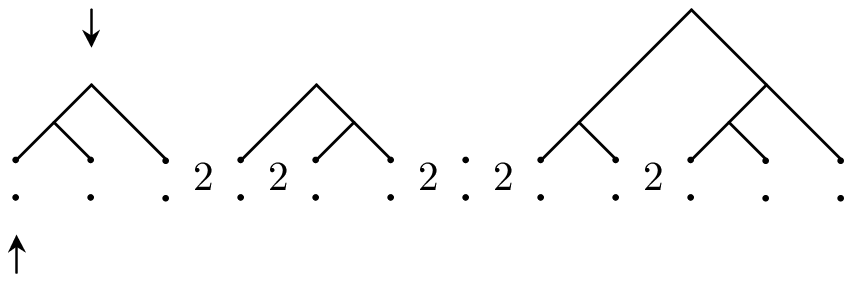}
\end{center}
Then $c(f)=8$ and $n(f)=5$, so $f$ has length $18$. The anti-normal form for~$f$~is:
\begin{equation*}
x_4x_5^2x_4x_2x_3x_1^2
\end{equation*}
Therefore, a minimum-length $\{x_0,x_1\}$-word for $f$ is:
\begin{equation*}
x_0^{-3}x_1x_0^{-1}x_1^2x_0x_1x_0^2x_1x_0^{-1}x_1x_0^2x_1^2\\[6pt]
\end{equation*}
\end{example}

Currently, our only algorithm to find the anti-normal form for a strongly
positive element involves drawing the forest diagram. It is interesting to
note that an entirely algebraic algorithm is available:

\begin{theorem}
Let $f\in F$ be strongly positive, and let $w$ be an expression for $f$
as a product of $\{x_1,x_2,\ldots\}$. Suppose we repeatedly apply operations
of the form:
\begin{equation*}
x_kx_n\longmapsto x_{n-1}x_k\parbox{0in}{\parbox{2in}{\qquad\emph{(}$k < n-1$\emph{)}}}
\end{equation*}
to $w$. Then we reach the anti-normal form for $f$ after at most
$\begin{pmatrix}c\left(f\right)\\2\end{pmatrix}$ steps.
\end{theorem}

\begin{proof}
Let $\mathcal{C}$ be the set of carets in the reduced forest diagram for $f$.  Suppose
that:
\begin{equation*}
w=x_{i_m}\cdots x_{i_2}x_{i_1}
\end{equation*}
Each generator $x_{i_k}$ appearing in $w$ corresponds to the construction of some caret
$c_k$ of the forest diagram for $f$.  Let $<$ denote the order in which these carets are created:
\begin{equation*}
c_1<c_2<\cdots<c_m
\end{equation*}
Now, the anti-normal form for $f$ is just another word for $f$ in the generators $\{x_1,x_2,\ldots\}$.
Let $<\!\!_{\scriptscriptstyle\text{AN}}$ denote the resulting order on $\mathcal{C}$.  Note that:
\begin{equation*}
c_k<\!\!_{\scriptscriptstyle\text{AN}}\:c_{k+1}\qquad\Longleftrightarrow\qquad
i_k - 1\leq i_{k+1}
\end{equation*}
Therefore, any operation of the form:
\begin{equation*}
x_{i_{k+1}}x_{i_k}\longmapsto
x_{i_k-1}x_{i_{k+1}}\parbox{0in}{\parbox{2in}{\qquad($i_{k+1}<i_k-1$)}}
\end{equation*}
reduces the number:
\begin{equation*}
\big|\left\{\left(c,c'\right):c<\!\!_{\scriptscriptstyle\text{AN}}\:c'\text{
but }c>c'\right\}\big|
\end{equation*}
by exactly one. When this number reaches zero, $f$ is in anti-normal form.

Finally, note that the number in question is bounded by
$\displaystyle\binom{\left|\mathcal{C}\right|}{2}$.
\end{proof}

\begin{example}
Let's find the length of the element:
\begin{equation*}
x_1\,x_3^3\,x_6\,x_7\,x_{10}
\end{equation*}
We put the word into anti-normal form:
\begin{equation*}
\begin{aligned}[b]
x_1\,&x_3^3\,x_6\,x_7\,x_{10}\\
= x_4\,&x_1\,x_3^3\,x_6\,x_7\\
= x_4\,&x_2^3\,x_5\,x_6\,x_1\\
= x_4\,&x_2\,x_3\,x_4\,x_2^2\,x_1
\end{aligned}
\parbox[b]{0in}{\parbox[b]{2in}{\qquad
$\begin{gathered}[b]
\\
\text{(}x_{10}\text{ moved left)}\\
\text{(}x_1\text{ moved right)}\\
\text{(}x_2^2\text{ moved right)}
\end{gathered}
$}}
\end{equation*}
Hence, the length is:
\begin{equation*}
\left(3+2+1+1+2+1+0\right)+7=17
\end{equation*}
\end{example}

\

\subsection{The Length Formula}

We now give the length formula for a general element of $F$. Afterwards, we
will give several examples to illustrate intuitively why the formula works.
We defer the proof to section 4.3.

Let $f\in F$, and let $\mathfrak{f}$ be its reduced forest diagram. We label
the spaces of each forest of $\mathfrak{f}$ as follows:
\begin{enumerate}
\item Label a space $\LL$ (for \emph{left}) if it exterior and to the
left of the pointer.
\item Label a space $\NN$ (for \emph{necessary}) if it lies immediately
to the left of some caret (and is not already labeled $\LL$).
\item Label a space $\RR$ (for \emph{right}) if it exterior and
to the right of the pointer (and not already labeled $\NN$).
\item Label a space $\II$ (for \emph{interior}) if it interior
(and not already labeled $\NN$).
\end{enumerate}
We assign a \emph{weight} to each space in the support of $\mathfrak{f}$
according to its labels:

\begin{center}
\begin{tabular}{c}
\\
$\text{top}$\\
$\text{label}$
\end{tabular}%
\begin{tabular}{|c|cccc|}
\multicolumn{5}{c}{$\text{bottom label}$}\\
\hline\rule{0ex}{2.5ex}
&$\LL$&$\NN$&$\RR$&$\II$\\
\hline
$\LL$&$2$&$1$&$1$&$1$\rule{0ex}{2.5ex}\\
$\NN$&$1$&$2$&$2$&$2$\\
$\RR$&$1$&$2$&$2$&$0$\\
$\II$&$1$&$2$&$0$&$0$\\
\hline
\end{tabular}
\end{center}

\begin{example}
Here are the labels and weights for a typical forest diagram:
\begin{center}
\includegraphics{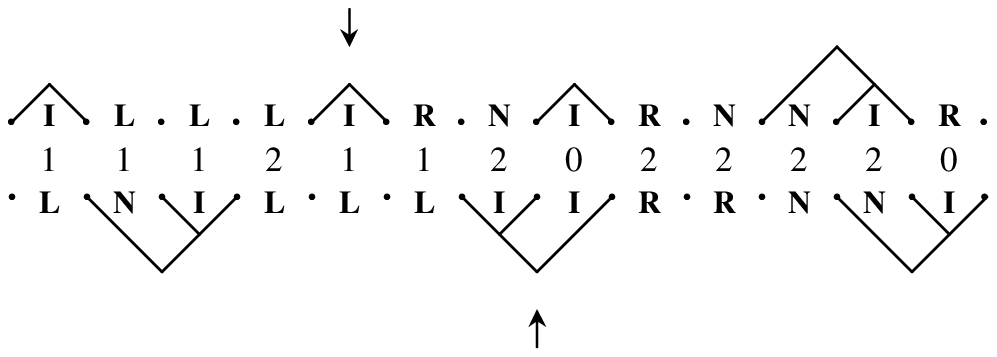}
\end{center}
\end{example}

\begin{theorem}[The Length Formula]
Let $f\in F$, and let $\mathfrak{f}$ be its reduced forest diagram. Then
the $\{x_0,x_1\}$-length of $f$ is:
\begin{equation*}
\ell(f)=\ell_0(f)+\ell_1(f)
\end{equation*}
where:
\begin{enumerate}
\item $\ell_0(f)$ is the sum of the weights of all spaces in the
support of\text{ }$\mathfrak{f}$, and
\item $\ell_1(f)$ is the total number of carets in $\mathfrak{f}$.
\end{enumerate}
\end{theorem}

\begin{remark}
Intuitively, the weight of a space is just the number of times it must
be crossed during the construction of $f$. Hence, there ought to exist
a minimum-length word for $f$ with $\ell_0(f)$ appearances of $x_0$ or $x_0^{-1}$
and $\ell_1(f)$ appearances of $x_1$ or $x_1^{-1}$. This will be established at the
end of the next section.
\end{remark}

\begin{example}
Let $f\in F$ be the element from example 4.1.9:
\begin{center}
\includegraphics{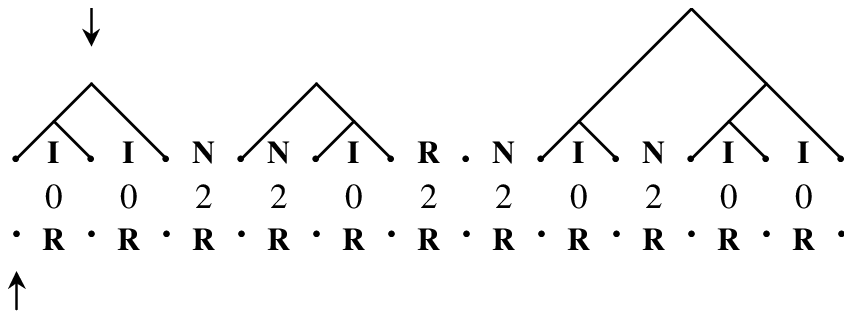}
\end{center}
Then $\ell_0(f)=10$ and $\ell_1(f)=8$, so $f$ has length $18$.\\[6pt]
\end{example}

In general, suppose $f\in F$ is strongly positive, and let $\mathfrak{f}$ be
its reduced forest diagram. Then every space of $\mathfrak{f}$ is labeled \pair{\NN}{\RR},
\pair{\RR}{\RR}, or \pair{\II}{\RR}.  Each \pair{\II}{\RR} space has weight $0$, and each
\pair{\NN}{\RR} or \pair{\RR}{\RR} space has weight $2$, so that:
\begin{equation*}
\ell_0(f)=2{\,}n(f)
\end{equation*}
and hence:
\begin{equation*}
\ell_0(f)+\ell_1(f)=2\,n(f)+c(f)
\end{equation*}
Therefore, the length formula of theorem 4.2.2 reduces to theorem 4.1.8
for strongly positive elements.

\begin{example}
Let $f$ be the inverse of the element from the previous example:
\begin{center}
\includegraphics{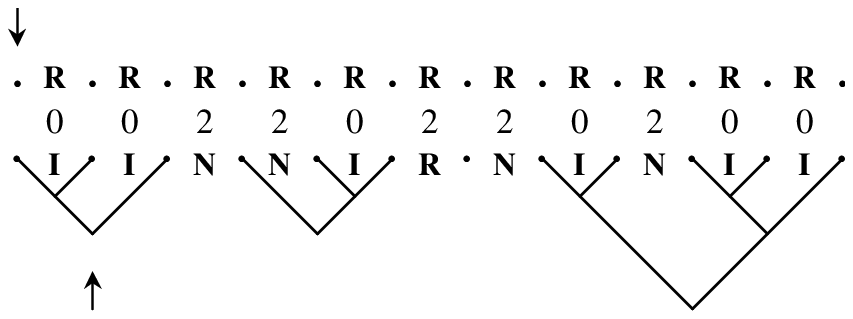}
\end{center}
Then $\ell_0(f)=10$ and $\ell_1(f)=8$, so $f$ has length $18$. One minimum-length word for $f$ is:
\begin{equation*}
x_1^{-2} x_0^{-2} x_1^{-1} x_0 x_1^{-1} x_0^{-2} x_1^{-1} x_0^{-1} x_1^{-2} x_0 x_1^{-1} x_0^3
\end{equation*}
Note that this word always creates the \emph{rightmost} possible caret first.
This is because the creation of a negative caret inserts a space above it,
and then moves the caret to the \emph{left} of this space. By creating
carets right to left, we avoid crossing over these newly created spaces.
\end{example}

\begin{example}
Let $f\in F$ be the element with forest diagram:
\begin{center}
\includegraphics{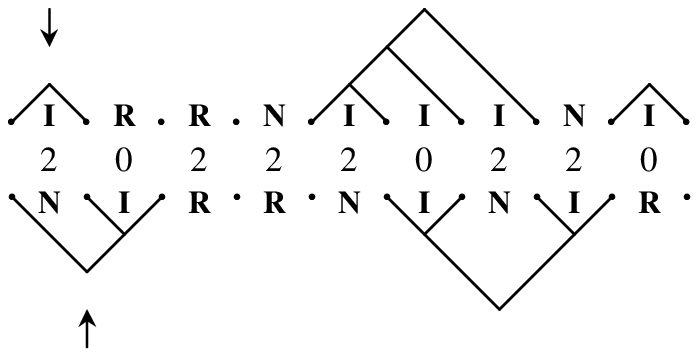}
\end{center}
Then $\ell_0(f)=12$ and $\ell_1(f)=10$, so $f$ has length $22$. One minimum-length word for $f$ is:
\begin{equation*}
x_1 x_0^{-1} x_1^{-1} x_0 x_1^{-1} x_0^{-3} x_1 x_0 x_1^3 x_0^{-1} x_1^{-1} x_0^{-1} x_1^{-1}
x_0 x_1^{-1} x_0^3\\[12pt]
\end{equation*}
\end{example}

In general, an element $f\in F$ is \emph{right-sided} if it lies in the subgroup
generated by $\{x_1,x_2,\ldots\}$. Equivalently, $f$ is right-sided if and only if
both pointers in the forest diagram for $f$ are at the left edge of the support.
Note then that every space of a right-sided element is labeled either $\NN$, $\RR$,
or $\II$. The weight table for such spaces is:
\begin{center}
\begin{tabular}{c}
\\
top\\
label
\end{tabular}%
\begin{tabular}{|c|ccc|}
\multicolumn{4}{c}{bottom label}\\
\hline\rule{0ex}{2.5ex}
&$\NN$&$\RR$&$\II$\\
\hline
$\NN$&$2$&$2$&$2$\rule{0ex}{2.5ex}\\
$\RR$&$2$&$2$&$0$\\
$\II$&$2$&$0$&$0$\\
\hline
\end{tabular}
\end{center}
Observe that a space has weight $2$ if and only if:
\begin{enumerate}
\item It is exterior on both the top and the bottom, or
\item It lies immediately to the left of some caret, on either the top or the bottom.
\end{enumerate}
This can be viewed as a generalization of the length formula for strongly positive elements.
Specifically, if $f$ is right-sided, then:
\begin{equation*}
\ell(f)=2\,n(f)+c(f)
\end{equation*}
where $n(f)$ is the number of spaces satisfying condition (1) or (2), and $c(f)$ is the number
of carets of $f$.

As with strongly positive elements, it is intuitively obvious that this is a
lower bound for the length. Unfortunately, we have not been able to find an
analogue of the ``anti-normal form'' argument to show that it is an upper
bound.

\pagebreak[4]
\begin{example}
Let $f\in F$ be the element with forest diagram:
\begin{center}
\includegraphics{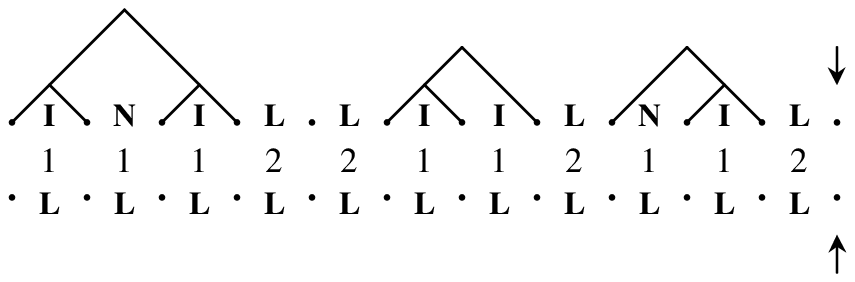}
\end{center}
Then $\ell_0(f)=15$ and $\ell_1(f)=7$, so $f$ has length $22$.

It is interesting to note that every interior space of $f$ has weight $1$:
for trees to the left of the pointer, one cannot avoid crossing interior
spaces at least once. Specifically, each caret is created from its
\emph{left} leaf, and we must move to this leaf somehow.

One minimum-length word for $f$ is
\begin{equation*}
x_0^4 x_1^2 x_0^{-2} x_1 x_0^{-3} x_1^2 x_0^{-3} x_1 x_0^{-1} x_1 x_0^{-2}
\end{equation*}
Note that this word creates carets \emph{right to left}.
\end{example}

\begin{example}
Let $f\in F$ be the element with forest diagram:
\begin{center}
\includegraphics{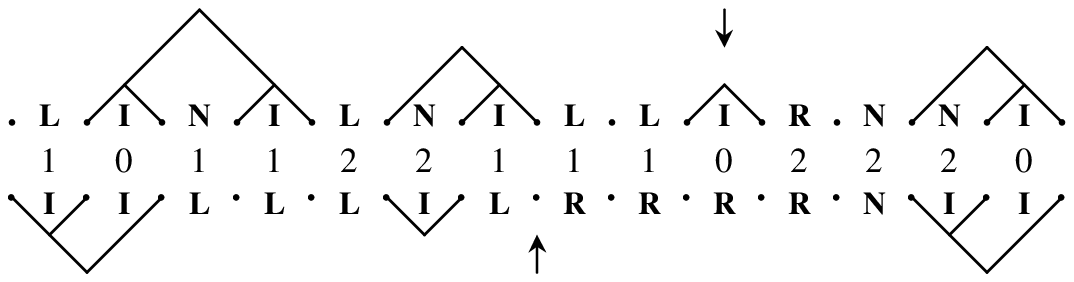}
\end{center}
Then $\ell_0(f)=16$ and $\ell_1(f)=13$, so $f$ has length $29$. One minimum-length word for $f$ is:
\begin{equation*}
x_0^{-2} x_1 x_0^{-1} x_1 x_0 x_1^{-2} x_0^2 x_1 x_0^3 x_1^2 x_0 x_1^{-2} x_0^{-1} x_1
x_0^{-2} x_1 x_0^{-1} x_1 x_0 x_1^{-1} x_0^{-1}
\end{equation*}
This is our first example with \pair{\LL}{\RR} pairs: note that they only need to be crossed once.
Also note how it affects the length to have bottom trees to the left of the pointer. In particular,
observe that the \pair{\NN}{\II} pair to the left of the pointers must crossed twice.
\end{example}

\

\subsection{The Proof of the Length Formula}

We prove the length formula using the same technique as Fordham \cite{Ford}:
\begin{theorem}
Let $G$ be a group with generating set $S$, and let $\ell\colon G\rightarrow\mathbb{N}$
be a function. Then $\ell$ is the length function for $G$ with respect to $S$ if and only if:
\begin{enumerate}
\item $\ell(e)=0$, where $e$ is the identity of $G$.
\item $\left|\ell(sg)-\ell(g)\right|\leq 1$ for all $g\in G$ and $s\in S$.
\item If $g\in G\setminus\{e\}$, there exists an $s\in S\cup S^{-1}$ such that $\ell(sg)<\ell(g)$.
\quad\qedsymbol
\end{enumerate}
\end{theorem}
\begin{proof}
Conditions (1) and (2) show that $\ell$ is a lower bound for the
length, and condition (3) shows that $\ell$ is an upper bound for the
length.
\end{proof}

Let $\ell$ denote the function defined on $F$ specified by Theorem 4.2.2.
Clearly $\ell$ satisfies condition (1). To show that $\ell$ satisfies
conditions (2) and (3), we need only gather information about how
left-multiplication by generators affects the function $\ell$.

\begin{terminology}
If $f\in F$, the \emph{current tree} of $f$ is the tree in forest diagram
indicated by the top pointer. The \emph{right space} of $f$ is the space
immediately to the right of the current tree, and the \emph{left space}
of $f$ is the space immediately to the left of the current tree.
\end{terminology}

\begin{proposition}
If $f\in F$, then $\ell(x_0f)=\ell(f)\pm 1$. Specifically, $\ell(x_0f)=\ell(f)-1$
unless one of the following conditions holds:
\begin{enumerate}
\item $x_0f$ has larger support than $f$.
\item The right space of $f$ has bottom label $\LL$, and left-multiplication
by $x_0$ does not remove this space from the support.
\item The right space of $f$ is labeled \pair{\RR}{\II}.
\end{enumerate}
\end{proposition}

\begin{proof}
Clearly $\ell_1(x_0f)=\ell_1(f)$. As for $\ell_0$, note that the only space whose
label changes is the right space of $f$.

\emph{Case 1}:\quad Suppose $x_0f$ has larger support than $f$. Then the right
space of $f$ is unlabeled, and has label \pair{\LL}{\RR} in $x_0f$. Hence
$\ell_0(x_0f)=\ell_0(f)+1$.

\emph{Case 2}:\quad Suppose $x_0f$ has smaller support than $f$. Then the right
space of $f$ has label \pair{\RR}{\LL}, but becomes unlabeled in $x_0f$. Hence
$\ell_0(x_0f)=\ell_0(f)-1$.

\emph{Case 3}:\quad Suppose $x_0f$ has the same support as $f$. Then the right
space of $f$ has top label $\NN$ or $\RR$, but top label $\LL$ in $x_0f$.
The relevant rows of the weight table are:
\begin{center}
\begin{tabular}{c}
\\
top\\
label
\end{tabular}%
\begin{tabular}{|c|cccc|}
\multicolumn{5}{c}{bottom label}\\
\hline\rule{0ex}{2.5ex}
&$\LL$&$\NN$&$\RR$&$\II$\\
\hline
$\LL$&$\boldsymbol{2}$&$1$&$1$\rule{0ex}{2.5ex}&$\boldsymbol{1}$\\
$\NN$&$\boldsymbol{1}$&$2$&$2$&$2$\\
$\RR$&$\boldsymbol{1}$&$2$&$2$&$\boldsymbol{0}$\\
\hline
\end{tabular}%
\end{center}
Each entry of the $\NN$ and $\RR$ rows differs from the corresponding entry
of the $\LL$ row by exactly one. In particular, moving from an $\RR$ or $\NN$
row to an $\LL$ row only decreases the weight when in the $\LL$ column or when
starting at~\pair{\RR}{\II}.
\end{proof}

\begin{corollary}
Let $f\in F$. Then $\ell(x_0^{-1}f)<\ell(f)$ if and only if one of the following conditions holds:
\begin{enumerate}
\item $x_0^{-1}f$ has smaller support than $f$.
\item The left space of $f$ has label \pair{\LL}{\LL}.
\item The left space of $f$ has label \pair{\LL}{\II}, and the current tree is trivial.
\end{enumerate}
\end{corollary}

\begin{proposition}
Let $f\in F$. If left-multiplying $f$ by $x_1$ cancels a caret from the
bottom forest, then $\ell(x_1f)=\ell(f)-1$.
\end{proposition}

\begin{proof}
Clearly $\ell_1(x_1f)=\ell_1(f)-1$. We must show that $\ell_0$ remains unchanged.

Note first that the right space of $f$ is destroyed. This space has label \pair{\RR}{\II},
and hence has weight $0$. Therefore, its destruction does not affect $\ell_0$.

The only other space affected is the left space of $f$. If this space is
not in the support of $f$, it remains unlabeled throughout. Otherwise,
observe that it must have top label $\LL$ in both $f$ and $x_1f$.
The relevant row of the weight table is:
\begin{equation*}
\begin{tabular}{|c|cccc|}
\hline\rule{0ex}{2.5ex}
&$\LL$&$\NN$&$\RR$&$\II$\\
\hline
$\LL$&$2$&$1$&$1$&$1$\rule{0ex}{2.5ex}\\
\hline
\end{tabular}%
\end{equation*}
In particular, the only important property of the bottom label is whether or
not it is an $\LL$. This property is unaffected by the deletion of the caret.
\end{proof}

\begin{proposition}
Let $f\in F$, and suppose that left-multiplying $f$ by $x_1$ creates a caret
in the top forest. Then $\ell(x_1f)=\ell(f)\pm 1$. Specifically, $\ell(x_1f)=\ell(f)-1$
if and only if the right space of $f$ has label \pair{\RR}{\RR}.
\end{proposition}

\begin{proof}
Clearly $\ell_1(x_1f)=\ell_1(f)+1$. As for $\ell_0$, observe that the only space whose
label could change is the right space of $f$.

\emph{Case 1}:\quad Suppose $x_1f$ has larger support than $f$. Then the right space of
$f$ is unlabeled, but has label \pair{\II}{\RR} in $x_1f$. This does not affect the value of $\ell_0$.

\emph{Case 2}:\quad Otherwise, note that the right space of $f$ has top label $\NN$ or
$\RR$. If the top label is an $\NN$, it remains and $\NN$ in $x_1f$. If it is an
$\RR$, then it changes to an $\II$. The relevant rows of the weight table are:
\begin{center}
\begin{tabular}{c}
\\
top\\
label
\end{tabular}%
\begin{tabular}{|c|cccc|}
\multicolumn{5}{c}{bottom label}\\
\hline\rule{0ex}{2.5ex}
&$\LL$&$\NN$&$\RR$&$\II$\\
\hline
$\RR$&$1$&$2$&$\boldsymbol{2}$&$0$\rule{0ex}{2.5ex}\\
$\II$&$1$&$2$&$\boldsymbol{0}$&$0$\\
\hline
\end{tabular}%
\end{center}
Observe that the weight decreases by two if the bottom label is an $\RR$,
and stays the same otherwise.
\end{proof}

We have now verified condition (2). Also, we have gathered enough
information to verify condition (3):

\begin{theorem}
Let $f\in F$ be a nonidentity element.
\begin{enumerate}
\item If current tree of $f$ is nontrivial, then either $\ell(x_1^{-1}f)<\ell(f)$,
or \mbox{$\ell(x_0f)<\ell(f)$}.
\item If left-multiplication by $x_1$ would remove a caret from the bottom tree, then 
$\ell(x_1f)<\ell(f)$.
\item Otherwise, either $\ell(x_0f)<\ell(f)$ or $\ell(x_0^{-1}f)<\ell(f)$.
\end{enumerate}
\end{theorem}

\begin{proof}

\

\emph{Statement 1}:\quad If $\ell(x_1^{-1}f)>\ell(f)$, then the right space of $x_1^{-1}f$
has type~\pair{\RR}{\RR}. The right space of $f$ therefore has type \pair{\RR\text{ or }\NN}{\RR\text{ or }\NN},
so that \mbox{$\ell(x_0f)<\ell(f)$}.

\emph{Statement 2}:\quad See proposition 5.3.5.

\emph{Statement 3}:\quad Suppose $\ell(x_0f)>\ell(f)$. There are three cases:

\emph{Case 1}:\quad The right space of $f$ is not in the support of $f$. Then the left space
of $f$ has label \pair{\LL}{\RR}, \pair{\LL}{\LL}, or \pair{\LL}{\II}. In all three cases,
$\ell(x_0^{-1}f)<\ell(f)$.

\emph{Case 2}:\quad The right space of $f$ has bottom label $\LL$, and
right-multiplication by $x_0$ does not remove this space from the support.
Then the left space of $f$ must have label \pair{\LL}{\LL} or \pair{\LL}{\II},
and hence $\ell(x_0^{-1}f)<\ell(f)$.

\emph{Case 3}:\quad The right space of $f$ has label \pair{\RR}{\II}. Then the tree immediately to the right
of the top pointer is trivial, and the bottom leaf under it is a right leaf. If the bottom leaf
under the top pointer were a left leaf, then left-multiplying $f$ by $x_1$ would cancel a caret.
Hence, it is also a right leaf, so the left space of $f$ has label \pair{\LL}{\II}. We conclude that
$\ell(x_0^{-1}f)<\ell(f).$
\end{proof}

\begin{corollary}
Let $f\in F$, and let $\mathfrak{f}$ be the reduced forest diagram for $f$. Then there exists
a minimum-length word $w$ for $f$ with the following properties:
\begin{enumerate}
\item Each instance of $x_1$ in $w$ creates a top caret of $\mathfrak{f}$.
\item Each instance of $x_1^{-1}$ in $w$ creates a bottom caret of $\mathfrak{f}$.
\end{enumerate}
In particular, $w$ has $\ell_1(f)$ instances of $x_1$ or $x_1^{-1}$, and $\ell_0(f)$
instances of $x_0$~or~$x_0^{-1}$.
\end{corollary}

\begin{proof}
By the previous theorem, it is always possible to travel from $f$ to the identity in such
a way that each left-multiplication by $x_1$ deletes a bottom caret and each left-multiplication
by $x_1^{-1}$ deletes a top caret.
\end{proof}

Of course, not every minimum-length word for $f$ is of the given form. We will
discuss this phenomenon in the next section.

\

\subsection{Minimum-Length Words}

In principle, the results from the last section specify an algorithm for
finding minimum-length words.  (Given an element, find a generator which
shortens it.  Repeat.)  In practice, though, no algorithm is necessary:
one can usually guess a minimum-length word by staring at the forest diagram.
Our goal in this section is to convey this intuition.

\begin{example}
Let $f$ be the element of $F$ with forest diagram:
\begin{center}
\includegraphics{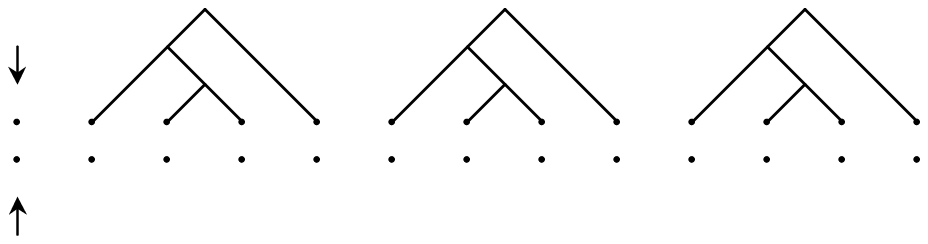}
\end{center}
Then there is exactly one minimum-length word for $f$, namely:
\begin{equation*}
x_0^{-3} u x_0 u x_0 u x_0
\end{equation*}
where $u = x_1^2 x_0^{-1} x_1 x_0$. Note that the trees of $f$ are constructed
from \emph{left to right}.

Similarly, $f^{-1}$ has forest diagram:
\begin{center}
\includegraphics{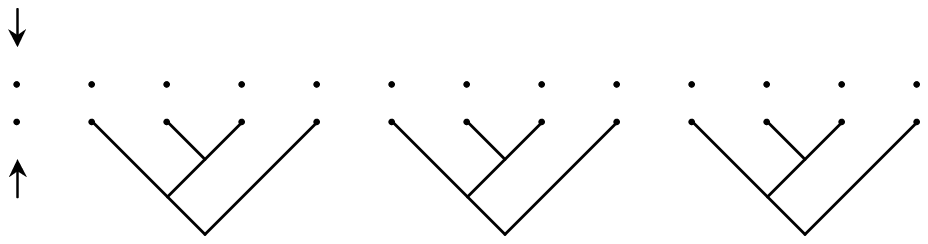}
\end{center}
and the only minimum-length word for $f^{-1}$ is:
\begin{equation*}
x_0^{-1} u^{-1} x_0^{-1} u^{-1} x_0^{-1} u^{-1} x_0^3
\end{equation*}
Note that the trees of $f^{-1}$ are constructed from \emph{right to left}.
\end{example}

\begin{example}
\quad Let $f$ be the element of $F$ with forest diagram:
\begin{center}
\includegraphics{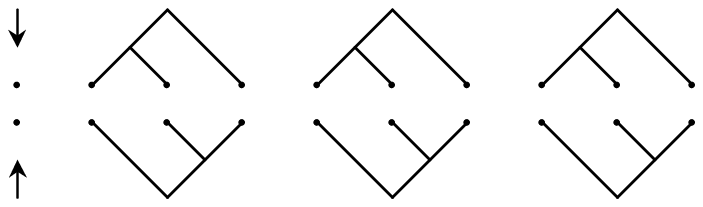}
\end{center}
There are precisely four minimum-length words for $f$:
\begin{equation*}
\begin{aligned}[t]
&x_0^{-4}\,v\,x_0\,v\,x_0\,v\,x_0\\
&x_0^{-1}\,v\,x_0^{-2}\,v\,x_0\,v\,x_0^2\\
&x_0^{-2}\,v\,x_0^{-1}\,v\,x_0^2\,v\,x_0\\
&x_0^{-1}\,v\,x_0^{-1}\,v\,x_0^{-1}\,v\,x_0^3
\end{aligned}
\end{equation*}
where $v = x_1^2 x_0^{-1} x_1^{-1} x_0 x_1^{-1}$. In particular, each of the
first two components can be constructed either when the pointer is moving
right, or later when the pointer is moving back left.
\end{example}

\begin{example}
Let $f$ be the element of $F$ with forest diagram:
\begin{center}
\includegraphics{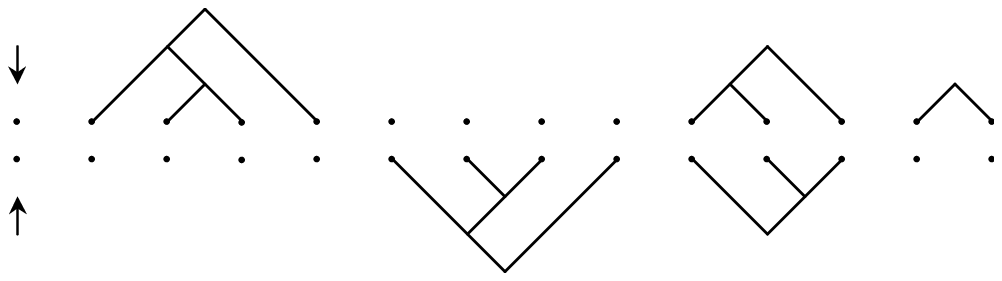}
\end{center}
There are precisely two minimum-length words for $f$:
\begin{equation*}
\begin{aligned}[t]
&x_0^{-2}\,u^{-1}\,x_0^{-2}\,x_1\,x_0\,v\,x_0^2\,u\,x_0\\
&x_0^{-2}\,u^{-1}\,x_0^{-1}\,v\,x_0^{-1}\,x_1\,x_0^3\,u\,x_0
\end{aligned}
\end{equation*}
where $u = x_1^2 x_0^{-1} x_1 x_0$ and $v = x_1^2 x_0^{-1} x_1^{-1} x_0 x_1^{-1}$. Note
that the first component must always be constructed on the journey right, and
the second component must always be constructed on the journey left. The
only choice lies with the construction of the third component: should it be
constructed when moving right, or should it be constructed while moving back
left?
\end{example}

In general, certain components act like ``top trees'' while others act like
``bottom trees'', while still others are ``balanced''. For example, the
forest diagram:
\begin{center}
\includegraphics{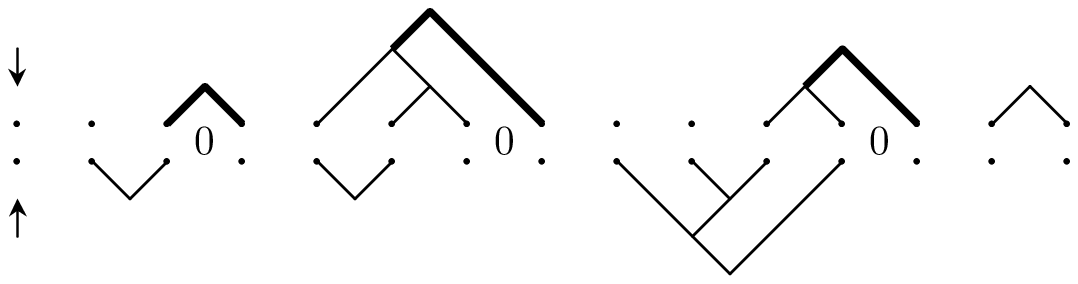}
\end{center}
must be constructed from left to right (so all the components act like ``top
trees''). The reason is that the three marked spaces each have weight $0$,
so that each of the three highlighted carets must be constructed \emph{before}
the pointer can move farther to the right. Essentially, the highlighted carets
are acting like \emph{bridges} over these spaces.

The idea of the ``bridge'' explains two phenomena we have already observed.
First, consider the following contrapositive of proposition 4.3.6:

\begin{proposition}
Let $f\in F$, and suppose that the top pointer of $f$ points at a
nontrivial tree. \ Then $\ell\left(x_1^{-1}f\right)<\ell(f)$
unless the resulting uncovered space has type \pair{\RR}{\RR}.\quad\qedsymbol
\end{proposition}

This proposition states conditions under which the destruction of a top caret
decreases the length of an element. Essentially, the content of the
proposition is that it makes sense to delete a top caret \emph{unless that caret is
functioning as a bridge}. (Note that the deletion of any of the bridges in
the example above would result in an \pair{\RR}{\RR} space.) It makes no sense to
delete a bridge, since the bridge is helping you access material further to the right.

Next, recall the statement of corollary 4.3.8:

\begin{corollary*}
Let $f\in F$, and let $\mathfrak{f}$ be the reduced forest diagram for
$f$. Then there exists a minimum-length word $w$ for $f$ with the following
properties:
\begin{enumerate}
\item Each instance of $x_1$ in $w$ creates a top caret of $\mathfrak{f}$.
\item Each instance of $x_1^{-1}$ in $w$ creates a bottom caret of
$\mathfrak{f}$.
\end{enumerate}
In particular, $w$ has $\ell_1(f)$ instances of $x_1$ or
$x_1^{-1}$, and $\ell_0(f)$ instances of $x_0$ or $x_0^{-1}$.
\end{corollary*}

When we originally stated this corollary, we mentioned that not every minimum-length word for $f$
is necessarily of the specified form. The reason is that it is sometimes reasonable to construct temporary
bridges while building an element.

\begin{example}
Let $f$ be the element of $F$ with forest diagram:
\begin{center}
\includegraphics{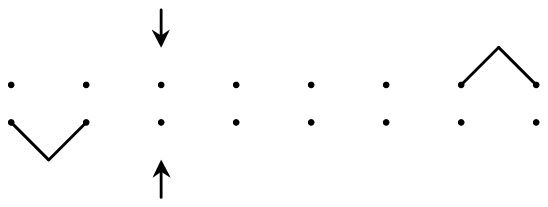}
\end{center}
Then one minimum-length word for $f$ is:
\begin{equation*}
x_0^2 x_1^{-1} x_0^{-5} x_1 x_0^4
\end{equation*}
This word corresponds to the instructions ``move right, create the top caret,
move left, create the bottom caret, and then move back to the origin''.
However, here is another minimum-length word for $f$:
\begin{equation*}
x_0^2 x_1^{-1} (x_0^{-1}x_1^{-3}x_0^{-1}) x_1 (x_0x_1^3)
\end{equation*}
In this word, the ``move right'' is accomplished by building three temporary
bridges:
\begin{center}
\includegraphics{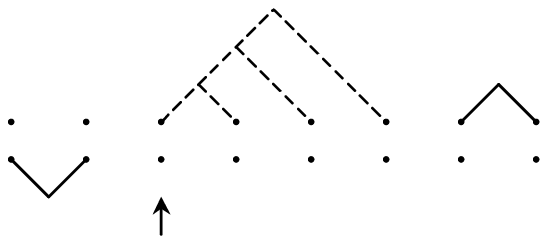}
\end{center}
These bridges are torn down during the ``move left''.

Finally, here is a third minimum-length word for $f$:
\begin{equation*}
x_1^{-3} x_0^2 x_1^{-1} x_0^{-2} x_1 (x_0x_1^3)
\end{equation*}
In this word, bridges are again built during the ``move right'', but they
aren't torn down until the very end of the construction.
\end{example}

We now turn our attention to a few examples with some more complicated
behavior.

\begin{example}
\quad Let $f$ be the element of $F$ with forest diagram:
\begin{center}
\includegraphics{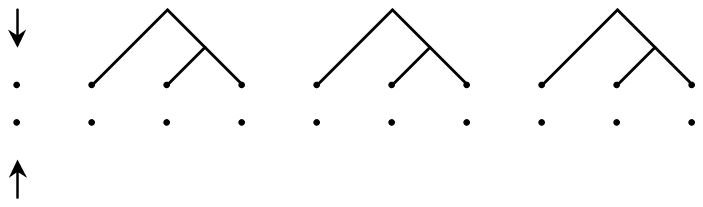}
\end{center}
There are four different minimum-length words for $f$:
\begin{equation*}
\begin{aligned}[t]
&x_0^{-3}x_1x_0^{-1}x_1x_0^2x_1x_0^{-1}x_1x_0^2x_1x_0^{-1}x_1x_0^2\\
&x_0^{-1}x_1x_0^{-3}x_1x_0^{-1}x_1x_0^2x_1x_0^{-1}x_1x_0^2x_1x_0^2\\
&x_0^{-2}x_1x_0^{-2}x_1x_0^{-1}x_1x_0^2x_1x_0^2x_1x_0^{-1}x_1x_0^2\\
&x_0^{-1}x_1x_0^{-2}x_1x_0^{-2}x_1x_0^{-1}x_1x_0^2x_1x_0^2x_1x_0^2
\end{aligned}
\end{equation*}
Note that each of the first two components may be either partially or fully
constructed during the move to the right. This occurs because the trees in
this example do not end with bridges. (Compare with example 4.4.1.)
\end{example}

\begin{example}
\quad Let $f$ be the element of $F$ with forest diagram:
\begin{center}
\includegraphics{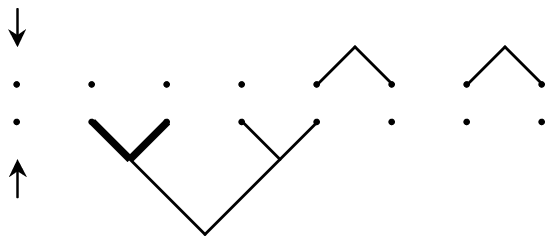}
\end{center}
There is exactly one minimum-length word for $f$:
\begin{equation*}
x_0^{-1} x_1^{-1} x_0^{-3} x_1 x_0 x_1 x_0 x_1^{-1} x_0 x_1^{-1} x_0
\end{equation*}
Note that the highlighted caret must be constructed \emph{last}, since the space it
spans should not be crossed. However, we must begin by partially constructing the
first component, because of the bridge on its right end.
\end{example}

\

\section{Applications}

This section contains various applications of forest diagrams and the length
formula.

\subsection{Dead Ends and Deep Pockets}

In \cite{ClTa1}, S. Cleary and J. Taback prove that $F$ has ``dead
ends'' but no ``deep pockets''. In this subsection, we show how forest
diagrams can be used to understand these results.

\begin{definition}
A \emph{dead end} is an element $f\in F$ such that
\mbox{$\ell(xf)<\ell(f)$} for all $x\in\left\{x_0, x_1, x_0^{-1}, x_1^{-1}\right\}$.
\end{definition}

\begin{example}
Consider the element $f\in F$ with forest diagram:
\begin{center}
\includegraphics{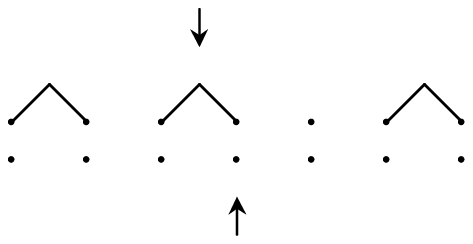}
\end{center}
Left-multiplying by $x_0^{-1}$ decreases the length since the left space of
$f$ is of type \pair{\LL}{\LL}. Left-multiplying by $x_0$ or $x_1$ decreases the
length since the right space of $f$ is of type \pair{\RR}{\RR}. Finally, left-multiplying
by $x_1^{-1}$ decreases the length since it deletes a top caret and the right
space of $x_1^{-1}f$ is not of type~\pair{\RR}{\RR}.
\end{example}

This example is typical:

\begin{proposition}
Let $f\in F$. Then $f$ is a dead end if and only if:
\begin{enumerate}
\item The current tree of $f$ is nontrivial,
\item The left space of $f$ has label \pair{\LL}{\LL},
\item The right space of $f$ has label \pair{\RR}{\RR}, and
\item The right space of $x_1^{-1}f$ does not have label \smash[t]{\pair{\RR}{\RR}}.
\end{enumerate}
\end{proposition}

\begin{proof}
The ``if'' direction is trivial. To prove the ``only if''
direction, assume that $f$ is a dead end. Then:
\begin{description}
\item[Condition (1)] follows from the fact that
$\ell\left(x_1^{-1}f\right)<\ell(f)$.
\item[Condition (2)] now follows from the fact that
$\ell\left(x_0^{-1}f\right)<\ell(f)$.
\item[Condition (3)] now follows from the fact that
$\ell(x_1f)<\ell(f)$.
\item[Condition (4)] now follows from the fact that
$\ell\left(x_1^{-1}f\right)<\ell(x_1f)$.\hfill\qedsymbol
\end{description}\renewcommand{\qedsymbol}{}
\end{proof}

Note that there are several ways to meet condition (4): the right space of
$x_1^{-1}f$ could be of type \pair{\RR}{\LL} (as in example 5.1.2), or it could
be of type \pair{\RR}{\II}:
\begin{center}
\includegraphics{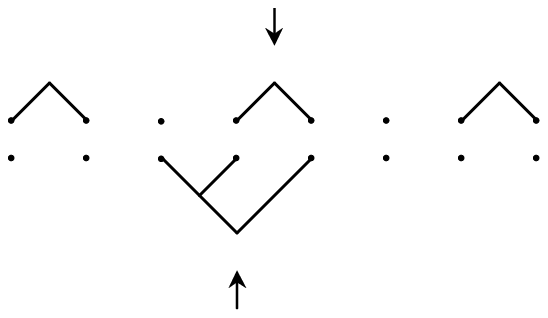}
\end{center}
or it could just have an $\NN$ on top:
\begin{center}
\includegraphics{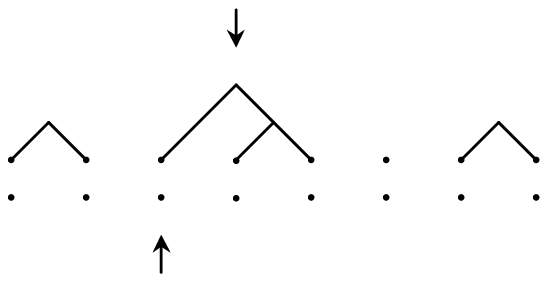}
\end{center}

\begin{definition}
Let $k\in\mathbb{N}$. A \emph{$k$-pocket} of $F$ is an element
$f\in F$ such that:
\begin{equation*}
\ell(s_1\cdots s_kf)\leq\ell(f)
\end{equation*}
for all $s_1,\ldots,s_k\in\left\{x_0,x_1,x_0^{-1},x_1^{-1},1\right\}$
\end{definition}

A 2-pocket in $F$ is just a dead end. S. Cleary and J. Taback demonstrated
that $F$ has no $k$-pockets for $k\geq 3$. We give an alternate proof:

\begin{proposition}
$F$ has no $k$-pockets for $k\geq 3$.
\end{proposition}

\begin{proof}
Let $f\in F$ be a dead-end element. Then the right space of $f$
has label \pair{\RR}{\RR}, so the tree to the right of the top pointer is trivial.
Therefore, repeatedly left-multiplying $x_0f$ by $x_1^{-1}$ will create
negative carets:
\begin{center}
\includegraphics{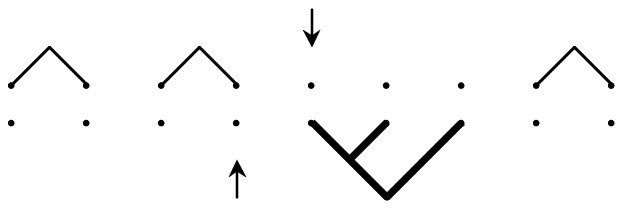}
\end{center}
In particular, $x_1^{-1}x_1^{-1}x_0f$ has length $\ell(f)+1$.
\end{proof}

\

\subsection{Growth}

We can use forest diagrams to calculate the growth function of the positive monoid with respect
to the $\{x_0,x_1\}$-generating set.  Burillo \cite{Bur} recently arrived at the same result using
tree diagrams and Fordham's length formula:

\begin{theorem}Let $P_n$ denote the number of positive elements of length $n$, and let:
\begin{equation*}
p(x) = \sum_{n=0}^\infty p_n x^n
\end{equation*}
Then:
\begin{equation*}
p(x) = \frac{1 - x^2}{1 - 2x - x^2 + x^3}
\end{equation*}
In particular, $p_n$ satisfies the recurrence relation:
\begin{equation*}
p_n = 2 p_{n-1} + p_{n-2} - p_{n-3}
\end{equation*}
for large $n$.
\end{theorem}

\begin{proof}
Let $P_n$ denote the set of all positive elements of length $n$.  Define four subsets of $P_n$ as follows:
\begin{enumerate}
\item $A_n = \{f\in P_n :$ the current tree of $f$ is trivial and is not the leftmost tree$\}$
\item $B_n = \{f\in P_n :$ the current tree of $f$ is nontrivial, but its right subtree is trivial$\}$
\item $C_n = \{f\in P_n :$ the current tree of $f$ is trivial and is the leftmost tree.$\}$
\item $D_n = \{f\in P_n :$ the current tree of $f$ is nontrivial, and so is its right subtree.$\}$
\end{enumerate}

Given an element of $A_n$, we can remove the current tree and move the pointer left, like this:
\begin{center}
\includegraphics{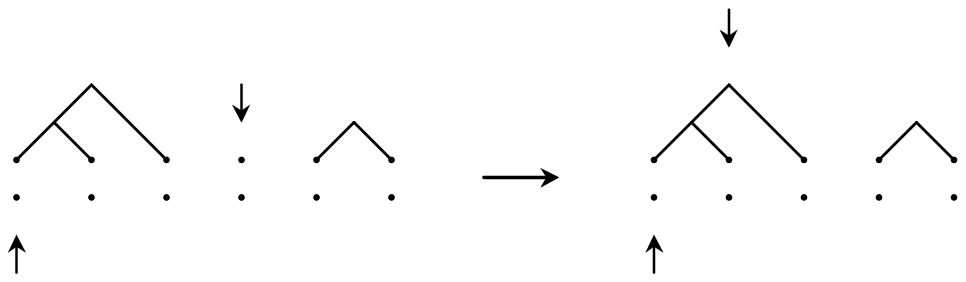}
\end{center}
This defines a bijection $A_n\rightarrow P_{n-1}$, so that:
\begin{equation*}
|A_n|=|P_{n-1}|
\end{equation*}

Given an element of $B_n$, we can remove the top caret together with the resulting trivial tree, like this:
\begin{center}
\includegraphics{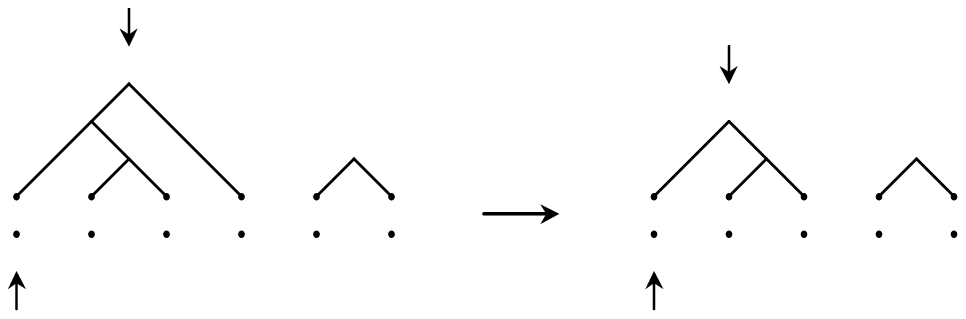}
\end{center}
This defines a bijection $B_n\rightarrow P_{n-1}$, so that:
\begin{equation*}
|B_n|=|P_{n-1}|
\end{equation*}

Given an element of $C_n$, we can move both the top and bottom arrows one space to the right, like this:
\begin{center}
\includegraphics{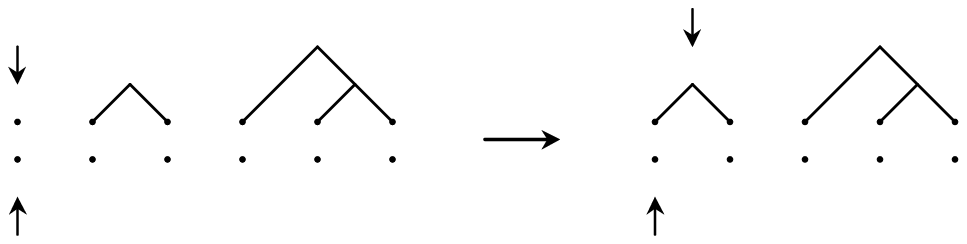}
\end{center}
This defines an injection $\varphi\colon C_n\rightarrow P_{n-2}$.  The image of $\varphi$ is all elements of
$P_{n-2}$ whose current tree is the first tree.

Finally, given an element of $D_n$, we can remove the top caret and move the pointer to the right subtree, like this:
\begin{center}
\includegraphics{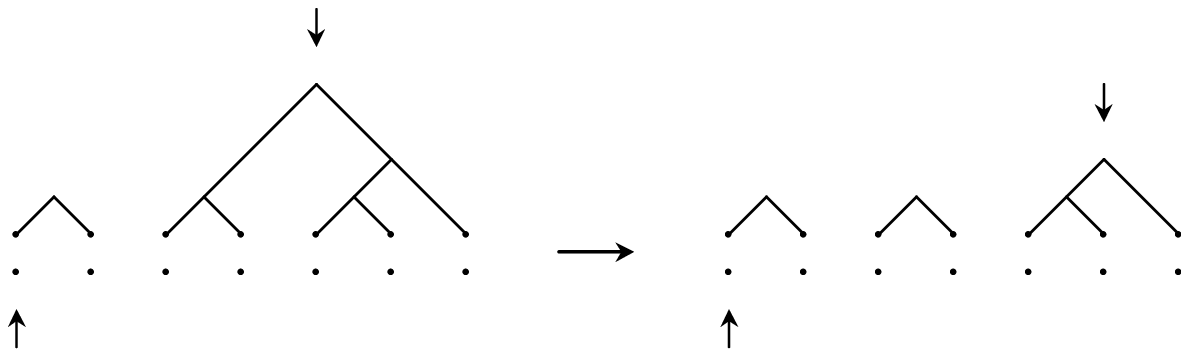}
\end{center}
This defines an injection $\psi\colon D_n\rightarrow P_{n-2}$.  The image of $\psi$ is all elements of
$P_{n-2}$ whose current tree is nontrivial, and is not the first tree.  In particular:
\begin{equation*}
(\text{im}\, \varphi) \cup (\text{im}\, \psi) = P_{n-2}-A_{n-2}
\end{equation*}
so that:
\begin{equation*}
|C_n|+|D_n|=|P_{n-2}|-|A_{n-2}|=|P_{n-2}|-|P_{n-3}|
\end{equation*}
This proves that $p_n$ satisfies the given recurrence relation for large $n$.  It is not much more work
to verify the given expression for $p(x)$.
\end{proof}

\subsection{The Isoperimetric Constant}

Let $G$ be a group with finite generating set $\Sigma$, and let $\Gamma$
denote the Cayley graph of $G$ with respect to $\Sigma$. If $S\subset G$,
define:
\begin{equation*}
\delta S=\left\{\text{edges in }\Gamma\text{ between }S\text{ and }S^c\right\}
\end{equation*}
The \emph{isoperimetric constant} of $G$ is defined as follows:
\begin{equation*}
\iota\left(G,\Sigma\right)=\inf\left\{\displaystyle\frac{\left|\delta
S\right|}{\left|S\right|}:S\subset G\text{ and }\left|S\right|<\infty\right\}
\end{equation*}
The group $G$ is amenable if and only if $\iota\left(G,\Sigma\right)=0$.

Guba \cite{Guba} recently proved that $\iota\left(F,\{x_0,x_1\}\right)\leq 1$. 
We have obtained a slightly better estimate:

\begin{proposition}
$\iota\left(F,\{x_0,x_1\}\right)\leq 1/2$.
\end{proposition}

\begin{proof}[Sketch of Proof]
Define the \emph{height} of a binary tree to be length of the
longest descending path starting at the root and ending at a leaf. Define the
\emph{width} of a binary forest to be the number of spaces in its support. For
each $n,k\in\mathbb{N}$, let $S_{n,k}$ denote all positive elements whose
forest diagram has width at most~$n$ and all of whose trees have height at most $k$.
One can show that:
\begin{equation*}
\lim_{k\rightarrow\infty}\lim_{n\rightarrow\infty}\frac{\left|\delta
S_{n,k}\right|}{\left|S_{n,k}\right|}=\frac{1}{2}
\parbox{0in}{\parbox{1.97in}{\raggedleft\qedsymbol}}
\end{equation*}\renewcommand{\qedsymbol}{}
\end{proof}

\

\subsection{Convexity}

A group $G$ is \emph{convex} (with respect to some generating set) if the $n$-ball $B^n(G)$ is a convex subset
of the Cayley graph of $G$ for each $n$.  Very few groups are convex, but Cannon \cite{Can} has introduced
the following weaker property:

\begin{definition}
A group $G$ is \emph{almost convex} (with respect to some generating set) if there exists an integer $L$
having the following property:  given any $x,y\in B^n(G)$ a distance two apart, there exists a path from
$x$ to $y$ in $B^n(G)$ of length at most $L$.
\end{definition}

The convexity of $F$ was first investigated by S. Cleary and J. Taback [ClTa2], who proved that $F$
is not almost convex with respect to $\{x_0,x_1\}$.  Recently, J. Belk and K. Bux \cite{BeBu} have applied forest
diagrams and the length formula to show that $F$ is \emph{maximally nonconvex}.  Specifically:

\begin{theorem}
For each $n\in\mathbb{N}$, let $l_n$ be the element of $F$ with forest diagram:
\begin{center}
\includegraphics{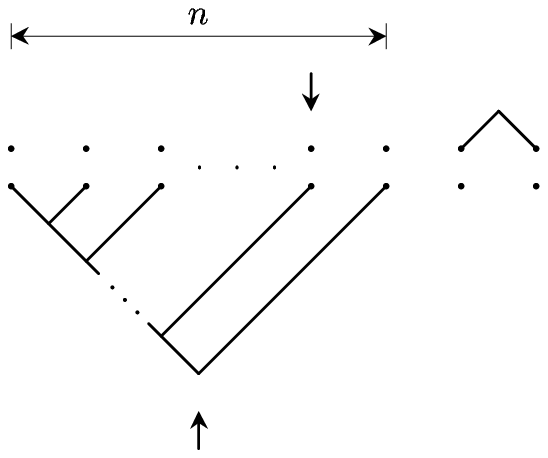}
\end{center}
and let $r_n = x_0^2 \, l_n$. Then $l_n$ and $r_n$ each have length $2n+2$,
and the shortest path from $l_n$ to $r_n$ inside the $(2n+2)$-ball has length $4n+4$.
\end{theorem}

\begin{proof}[Sketch of Proof]:
Since the right space of $l_n$ has label \pair{\RR}{\II}, $x_0 l_n$ has greater length than $l_n$:
\begin{equation*}
\includegraphics{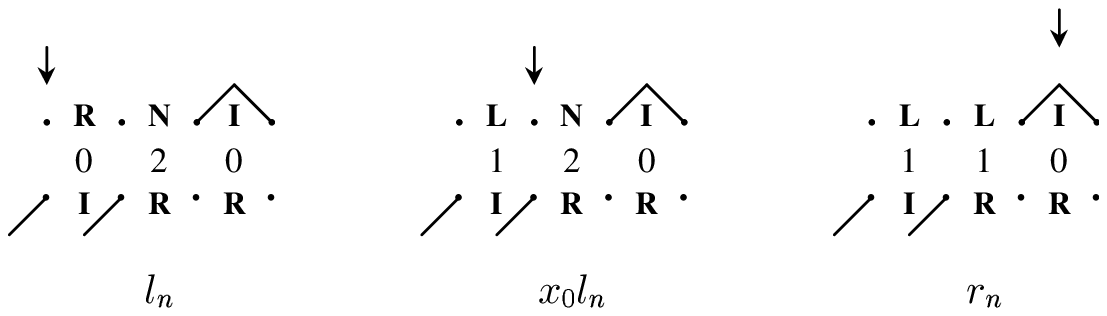}
\end{equation*}
In particular, the path:
\begin{equation*}
l_n\;\text{---}\;x_0l_n\;\text{---}\;r_n
\end{equation*}
does not remain within the $(2n+2)$-ball.

Intuitively, if one wants to get from $l_n$ to $r_n$ while remaining inside the $(2n+2)$-ball,
one must begin by moving all the way to the left and removing the accessible bottom caret.
Taking this idea further, we might guess that the following path of length $4n+4$ is minimal:
\begin{enumerate}
\item Move left $n-1$ spaces, and delete the leftmost bottom caret.
\item Move right $n$ spaces, and delete the top caret.
\item Move left $n$ spaces, and re-create the leftmost bottom caret.
\item Move right $n+1$ spaces, and re-create the top caret.
\end{enumerate}
This is in fact the case (see \cite{BeBu}).
\end{proof}

\

\end{document}